\documentclass[12pt]{article}

\usepackage{graphicx}
\graphicspath{{dir_pm_fig/}}
\usepackage{epstopdf}
\usepackage{amsopn}
\usepackage{amsmath}
\usepackage{amssymb}
\usepackage{amsfonts}
\usepackage[margin=1.5cm]{geometry}
\usepackage{siunitx}
\newcommand{\imunit}{\mathrm{i}}
\newcommand{\euler}{\mathrm{e}}
\newcommand{\Real}{\mathbb{R}}
\newcommand{\Complex}{\mathbb{C}}
\newcommand{\Integer}{\mathbb{Z}}
\newcommand{\Leb}{L}

\newcommand\given[1][]{\:#1\vert\:}
\newcommand{\vect}[1]{\boldsymbol{#1}}
\newcommand{\fourier}[1]{\hat{#1}}
\newcommand{\bessel}[1]{\lowercase{#1}}
\newcommand{\wfourier}[1]{\widehat{#1}}
\newcommand{\rotation}{{\cal R}}
\newcommand{\translation}{{\cal T}}

\newcommand{\transpose}[1]{{#1}^{\intercal}}
\newcommand{\vx}{{\vect{x}}}
\newcommand{\vk}{{\vect{k}}}
\newcommand{\vd}{{\vect{\delta}}}
\newcommand{\vu}{{\vect{u}}}
\newcommand{\vv}{{\vect{v}}}
\newcommand{\vg}{{\vect{g}}}
\newcommand{\cX}{{\cal X}}

\newcommand{\cN}{{\cal N}}
\newcommand{\fA}{{\fourier{A}}}
\newcommand{\fB}{{\fourier{B}}}
\newcommand{\fsigma}{{\fourier{\sigma}}}
\newcommand{\fS}{{\fourier{S}}}
\newcommand{\bA}{{\bessel{A}}}
\newcommand{\bB}{{\bessel{B}}}

\newcommand{\kmax}{{K}}
\newcommand{\qmax}{{Q}}
\newcommand{\rmax}{{R}}
\newcommand{\lmax}{{L}}
\newcommand{\mmax}{{M}}
\newcommand{\dx}{{\Delta x}}     % the pixel size in spatial-units
\newcommand{\dk}{{\Delta k}}     % the pixel size in frequency-units
\newcommand{\dpsi}{{\Delta \psi}}     % the pixel size in spatial-units
\newcommand{\Ccost}{C}
\newcommand{\Ckern}{\boldsymbol{C}}
\newcommand{\Dcost}{D}
\newcommand{\Dkern}{\boldsymbol{D}}

\newcommand{\argmax}{\operatorname*{arg\,max}}

\newcommand{\bigO}{\mathcal{O}}

\newcommand{\nimage}{N_{A}}
\newcommand{\ntarget}{N_{B}}

\newcommand{\npixel}{N}
\newcommand{\nrank}{H}
\newcommand{\signal}{{\text{\tiny signal}}}
\newcommand{\noise}{{\text{\tiny noise}}}
\newcommand{\optimal}{{\text{\tiny optimal}}}
\newcommand{\estimated}{{\text{\tiny estim}}}

\newcommand{\hk}{\hat{\boldsymbol{k}}}
\newcommand{\kpolara}{\theta}
\newcommand{\kazimub}{\phi}

\newcommand{\epolara}{\beta}
\newcommand{\eazimub}{\alpha}
\newcommand{\egammaz}{\gamma}

\begin{document}

% short title is in [..], full in {..} :
\title{Radial-recombination for rigid rotational alignment of images and volumes}
\author{Aaditya Rangan$^{1}$}

%\address{$^1$ Courant Institute, New York University, New York, NY and Center for Computational Mathematics, Flatiron Institute, New York, NY} \ead{rangan@cims.nyu.edu}

\begin{abstract}
% 276 words. ;
A common task in single particle electron cryomicroscopy (cryo-EM) is the rigid alignment of images and/or volumes.
In the context of images, a rigid alignment involves estimating the inner-product between one image of $\npixel\times\npixel$ pixels and another image that has been translated by some displacement and rotated by some angle $\gamma$.
In many situations the number of rotations $\gamma$ considered is large (e.g., $\bigO(\npixel)$), while the number of translations considered is much smaller (e.g., $\bigO(1)$).
In these scenarios a naive algorithm requires $\bigO(\npixel^{3})$ operations to calculate the array of inner-products for each image-pair.
This computation can be accelerated by using a fourier--bessel basis and the fast-fourier-transform (FFT), requiring only $\bigO(\npixel^{2})$ operations per image-pair.
We propose a simple data-driven compression algorithm to further accelerate this computation, which we refer to as the `radial-SVD'.
Our approach involves linearly-recombining the different rings of the original images (expressed in polar-coordinates), taking advantage of the singular-value-decomposition (SVD) to choose a low-rank combination which both compresses the images and optimizes a certain measure of angular discriminability.
When aligning multiple images to multiple targets, the complexity of our approach is $\bigO(\npixel(\log(\npixel)+\nrank))$ per image-pair, where $\nrank$ is the rank of the SVD used in the compression above.
A very similar strategy can be used to accelerate volume-alignment, using a spherical-harmonic based compression, which we'll refer to as a `degree-SVD'.
The advantage gained by this approach depends on the ratio between $\nrank$ and $\npixel$; the smaller $\nrank$ is the better.
In many applications $\nrank$ can be quite a bit smaller than $\npixel$ while still maintaining accuracy.
We present numerical results in a cryo-EM application demonstrating that the radial- and degree-SVD can help save a factor of $5$--$10$ for both image- and volume-alignment.
\end{abstract}

%\vspace{2pc}
%\noindent{\it Keywords}: single-particle electron cryomicroscopy, rigid image alignment, Fourier--Bessel basis, rigid volume alignment, Spherical--Harmonic basis, singular value decomposition

\maketitle

\section{Introduction}
\label{sec_Introduction}

Rigid alignment of images and volumes is a ubiquitous task that arises in many computer vision problems, as well as in the analysis of biomedical data \cite{CapekPousek2003,Szeliski2006,Shatsky2009}.
In this paper we focus on the application of single particle electron cryomicroscopy (cryo-EM) (see \cite{Cheng2015,Nogales2015,Elmlund2015,Murata2017,Sigworth2016} for an overview).

Within this application there are two basic tasks, involving the comparison of either: (i) a pair of 2-dimensional images, or (ii) a pair of 3-dimensional volumes.
During cryo-EM molecular-reconstruction the same images or volumes are typically compared multiple times against different targets \cite{Goncharov1987,vanHeel1987,Goncharov1988,Jonic2005,Grigorieff2007,EMAN2,Yang2008,Singer2009,Singer2011,Scheres2012,Shkolnisky2012,Lyumkis2013,Wang2013,Grigorieff2016,Punjani2017,PunjaniBrubaker2017}.
The images involved are typically sampled on a uniform square $\npixel\times\npixel$ grid, while the volumes are measured on a uniform cubic $\npixel\times\npixel\times\npixel$ grid.
The measure of similarity is typically the inner-product between the two objects, which is equivalent (up to scaling) to their correlation.
In each case the goal is to find a rigid transformation -- a rotation composed with a translation -- that maximizes the inner-product between one object and the other.

In the cryo-EM setting the images and volumes involved are usually approximately centered to begin with.
Thus, while the set of transformations to be considered still includes all rotations, the translations can be limited to small translations with a bounded magnitude.
In this setting the translational degrees-of-freedom can be treated efficiently \cite{RSAB20}, and the problem reduces to a sub-problem where only rigid-rotations need be considered.

Within the context of images, the main challenge then becomes to solve the rotational-alignment problem quickly.
In some applications only the optimal rotation is sought, while in other applications the entire landscape of inner-products is required, as sampled over rotation-angles $\in [0,2\pi)$ \cite{Scheres2012,Punjani2017}.
As described in \cite{Joyeux2002,Barnett2017}, the computational complexity of the latter task is no greater than the former, and can be performed efficiently by representing the images in a fourier-bessel basis and using the 1-dimensional fast-fourier-transform (FFT).
For volumes, a very similar strategy can be applied: representing the volumes in a spherical-harmonic basis and using the 2-dimensional FFT \cite{Kostelec03fftson}.
In all of these cases the computational efficiency of the calculation depends on representing the objects -- either images or volumes -- in a basis that explicitly records the radial-direction (i.e., either image-rings or volume-shells, respectively).

In this paper we make a very simple observation: not all of the different radii are equally useful.
Some are relevant for discriminating between the objects, while others carry little to no information.
Motivated by this observation, we propose linearly recombining the image-rings (or volume-shells, resp.) to reduce the number of degrees-of-freedom, while maintaining information relevant for alignment.
This can easily be done by constructing an (application-dependent) quadratic objective-function which measures the quality of any given radial-recombination.
This objective-function can then be optimized using the singular-value-decomposition (SVD), revealing radial-combinations that are useful for alignment.
As we demonstrate below, this `radial-SVD' can be used to reduce the operation-count of image- and volume-alignment, often by a factor of $5$--$10$ or more when aligning multiple images or volumes to multiple targets.

This paper is structured as follows.
We first review some mathematical preliminaries, introducing our notation and referencing the standard strategies for the rotational-alignment of images.
Then we describe our proposed radial-SVD in the context of images, and present an example using cryo-EM data.
We then generalize our approach to volumes, using the degree-SVD as well as the radial-SVD.
We conclude by pointing out some of the obvious applications for the radial-SVD, as well as some straightforward generalizations.

\section{Mathematical Preliminaries}
\label{sec_Mathematical_Preliminaries}

In this section we introduce the notation used in the rest of the manuscript, describe the objects we are dealing with in the context of cryo-EM, and review the standard strategies for calculating inner-products.
We first discuss images, and then summarize the generalization to volumes in section \ref{sec_Volume_alignment}.
Many of the conventions we establish for the former will carry over to the latter.
When possible, we will use the same notation as in \cite{RSAB20}.

\subsection{Image notation}
\label{sec_Image_notation}

We use $\vx,\vk \in\Real^{2}$ to represent spatial position and frequency, respectively.
In polar-coordinates these vectors are represented as:
%In polar-coordinates the vector $\vk$ is represented as:
\begin{eqnarray}
\vx &=& (x \cos \theta, x \sin \theta) \\
\vk &=& (k \cos \psi, k \sin \psi) \text{.}
\end{eqnarray}

The fourier transform of a two-dimensional function $A \in \Leb^{2}(\Real^{2})$ is defined as
\begin{equation}
  \fA(\vk) := \iint_{\Real^{2}} A(\vx) \euler^{-\imunit \vk \cdot \vx} \mathop{d\vx}~.
  \label{eq_Ahat}
\end{equation}
We recover $A$ from $\fA$ using the inverse fourier transform:
\begin{equation}
A(\vx) = \frac{1}{(2\pi)^{2}} \iint_{\Real^{2}} \fA(\vk) \euler^{+\imunit \vk \cdot \vx} \mathop{d\vk} \text{.}
\end{equation}
The inner-product between two functions $A, B \in \Leb^{2}(\Real^{2})$ is written as
\begin{equation}
  \langle A, B\rangle = \iint_{\Real^{2}} A(\vx)^{\dagger} B(\vx) \mathop{d\vx}
  ~,
\end{equation}
where $z^{\dagger}$ is the complex conjugate of $z \in \Complex$.
We will also use Plancherel's theorem \cite{bracewell},
\begin{equation}
  \langle A, B\rangle = \frac{1}{(2\pi)^{2}} \langle\fA,\fB\rangle
  ~, \qquad \forall A, B \in \Leb^{2}(\Real^{2})~.
  \label{eq_plancherel}
\end{equation}

We represent any given image as a function $A \in \Leb^{2}(\Real^{2})$, with values corresponding to the image intensity at each location.
As a consequence of Plancherel's theorem, any inner-product between $A$ and $B$ can be calculated equally well in either real- or frequency-space.
Because of recent image-alignment tools developed in the context of cryo-EM molecular-reconstruction \cite{Barnett2017,RSAB20}, we will typically refer to images in frequency-space (i.e., $\fA$, rather than $A$).
Nevertheless, all of the concepts we develop can be applied just as well in real-space \cite{Yang2008,Zhao2014,Sigworth2016,Barnett2017,RSAB20}.

Abusing notation, we'll refer to $A(\vk)$ and $\fA(\vk)$ in polar-coordinates as:
\begin{eqnarray}
A(x, \theta) &:=& A(\vx) = A(x \cos \theta, x \sin \theta) \\
\fA(k, \psi) &:=& \fA(\vk) = \fA(k \cos \psi, k \sin \psi) \text{.}
\end{eqnarray}

With this notation each $\fA(k,\psi)$ for fixed $k$ and $\psi\in[0,2\pi)$ corresponds to a `ring' in frequency-space with radius $k$.
%Later on in section \ref{sec_Radial_svd} we will consider linear-combinations of these image-rings.

\subsection{Image inner-products in a continuous setting: the fourier-bessel basis}
\label{sec_Image_innerproducts_in_a_continuous_setting}

Using the notation above, a rotation $\rotation_{\gamma}$ by angle $\gamma$ can be represented as:
\begin{equation}
\label{eq_rotation_definition}
\rotation_{\gamma} A(x, \theta) := A(x, \theta-\gamma) \text{.}
\end{equation}
Since rotation commutes with the fourier-transform, we have:
\begin{equation}
\wfourier{\rotation_\gamma A}(k, \psi) = \rotation_\gamma \circ \fA (k, \psi) = \fA(k, \psi-\gamma) \text{.}
\end{equation}
In this manner, a rotation of any image by $+\gamma$ can be represented as an angular-shift of each image-ring by $\psi\rightarrow\psi-\gamma$.

The inner-product between an image $A$ and a rotated-version of image $B$ is denoted by:
\begin{eqnarray}
\cX(\gamma ; A,B) &:=& \langle A , \rotation_{\gamma} B \rangle \text{,}
\end{eqnarray}
and is equivalent (up to a prefactor, which we ignore) to:
\begin{eqnarray}
%\cX(\gamma ; \fA,\fB) & =& \iint_{\Omega_{1}} A(x,\theta)^{\dagger} B(x,\theta-\gamma) xdxd\theta \\
%                 & =& \iint_{\Omega_{K}} \fA(k,\psi)^{\dagger} \fB(k,\psi-\gamma) kdkd\psi \text{.}
%\cX(\gamma ; \fA,\fB) & =& \iint A(x,\theta)^{\dagger} B(x,\theta-\gamma) xdxd\theta \\
\cX(\gamma ; \fA,\fB) & =& \langle \fA , \rotation_{\gamma}\fB \rangle \\
                 & =& \iint \fA(k,\psi)^{\dagger} \fB(k,\psi-\gamma) kdkd\psi \text{.}
\end{eqnarray}

It is computationally expensive to calculate many such inner-products $\cX(\gamma)$ (across multiple $\gamma$) directly from the values of $\fA$ and $\fB$, because each unique $\gamma$ corresponds to a different shift in the angular-component of each image-ring of $\fB$.
Put another way, the calculation of all $\cX(\gamma)$ for $\gamma\in[0,2\pi)$ amounts to a $\psi$-convolution of the image-rings $\fA(k,\psi)$ and $\fB(k,\psi)$.
To calculate such a convolution efficiently, we represent the image-rings in a fourier-bessel basis \cite{Zhao2014,Zhao2016,Barnett2017,RSAB20}.
This basis transforms the action of the rotation-operator into elementwise multiplication, and transforms the convolution above into an elementwise product.

To define the fourier-bessel-coefficients of an image we recall that, for each fixed $k$, the image-ring $\fA(k, \psi)$ is a $2\pi$-periodic function of $\psi$.
Thus, we can represent each image-ring $\fA(k, \psi)$ as a fourier-series in $\psi$, obtaining
\begin{equation}
\fA(k, \psi) = \sum_{q=-\infty}^{+\infty} \bA(k;q) \euler^{+\imunit q \psi} \text{,}
\end{equation}
for $q\in\Integer$.
The fourier-bessel-coefficients $\bA(k;q)$ of the image-ring $\fA(k,\psi)$ are given by
\begin{equation}
\bA(k;q) = \frac{1}{2\pi} \int_{0}^{2\pi} \fA(k, \psi) \euler^{-\imunit q \psi} \mathop{d\psi}~.
\label{eq_aqk}
\end{equation}
These coefficients can be represented in a more traditional fashion by recalling that the bessel-function $J_{q}(kx)$ can be written as:
\begin{eqnarray}
J_{q}(kx) & = & \frac{1}{2\pi}\int_{0}^{2\pi} \euler^{\imunit kx\sin(\psi)-\imunit q \psi} \mathop{d\psi}~ \\
 & = & \frac{1}{2\pi}\int_{0}^{2\pi} \euler^{-\imunit kx\cos(\psi+\pi/2)-\imunit q \psi} \mathop{d\psi}~ ,
\end{eqnarray}
which, when combined with the definition of the fourier-transform, immediately implies that:
\begin{eqnarray}
\bA(k;q) & = & \iint A(x, \theta) \frac{1}{2\pi} \int_{0}^{2\pi} \euler^{-\imunit kx\cos(\psi-\theta) - \imunit q \psi} \mathop{d\psi} xdx\mathop{d\theta}~ \\
 & = & \iint A(x, \theta) \euler^{\imunit(\theta+\pi/2)} \frac{1}{2\pi} \int_{0}^{2\pi} \euler^{-\imunit kx\cos(\psi+\pi/2) - \imunit q \psi} \mathop{d\psi} xdx\mathop{d\theta}~\\
 & = & \iint A(x, \theta) \euler^{\imunit(\theta+\pi/2)}J_{q}(kx) xdx\mathop{d\theta}~,
\end{eqnarray}
which is the inner-product between the original image (in real-space) and a `fourier-bessel' function.

The rotation $\rotation_{\gamma}$ can now be represented as:
\begin{eqnarray}
\rotation_{\gamma}\fA(k, \psi) & = & \fA(k,\psi-\gamma) \\
                      & = & \sum_{q=-\infty}^{+\infty} \bA(k;q) \euler^{+\imunit q (\psi-\gamma)} \\
                      & = & \sum_{q=-\infty}^{+\infty} \bA(k;q)\cdot\euler^{-\imunit q\gamma} \euler^{+\imunit q \psi} \text{,}
\label{eq_fourier_bessel_rotation}
\end{eqnarray}
such that the fourier-bessel-coefficients of the rotated image-ring $\rotation_{\gamma}\circ \fA(k,\cdot)$ are given by the original fourier-bessel-coefficients $\bA(k;q)$, each multiplied by the phase-factor $\euler^{-\imunit q\gamma}$.
Note that \eqref{eq_fourier_bessel_rotation} naturally allows for rotation by any $\gamma$, even values of $\gamma$ that may not lie on the polar grid used to discretize the images.

Using Plancherel's theorem in 1-dimension, we see that, for any pair of image-rings, 
\begin{equation}
\int_{0}^{2\pi} \fA(k,\psi)^{\dagger} \fB(k,\psi) d\psi = 2\pi \sum_{q=-\infty}^{+\infty} \bA(k,q)^{\dagger}\bB(k,q) \text{.}
\end{equation}
Thus, the inner-product $\cX(\gamma)$ can be written in terms of fourier-bessel-coefficients of $\fA$ and $\fB$:
\begin{eqnarray}
\cX(\gamma ; \fA,\fB) & = & \langle \fA , \rotation_{\gamma}\fB \rangle \\
%            & = & \int_{\Omega_{\kmax}} \sum_{q=-\infty}^{+\infty} \bA(k;q)^{\dagger}\bB(k;q) \euler^{-\imunit q \gamma} kdk \\
            & = & 2\pi \int \sum_{q=-\infty}^{+\infty} \bA(k;q)^{\dagger}\bB(k;q) \euler^{-\imunit q \gamma} kdk \\
            & = & 2\pi \sum_{q=-\infty}^{+\infty} \euler^{-\imunit q \gamma} \cdot \left[ \int  \bA(k;q)^{\dagger}\bB(k;q) kdk \right] \text{.}
\end{eqnarray}
This last expression can be interpreted as a relationship between the desired inner-products $\cX(\gamma)$ and the fourier transform of the term in brackets on the right-hand-side. That is:
\begin{equation}
\fourier{\cX}(q ; \bA,\bB) = 2\pi \left[ \int  \bA(k;q)^{\dagger}\bB(k;q) kdk \right] \text{.}
\end{equation}

In a typical discretization scheme (see below) the images $\fA$ and $\fB$ each require the storage of $\bigO(\npixel^{2})$ values, as do the fourier-bessel-representations $\bA$ and $\bB$.
The calculation of the array $\fourier{\cX}(q)$ involves $\bigO(\npixel^{2})$ operations and $\bigO(\npixel)$ storage.
Once $\fourier{\cX}$ is calculated, the inner-products $\cX(\gamma)$ can be recovered on a uniform grid of $\bigO(\npixel)$ angles $\gamma\in[0,2\pi)$ with an additional $\bigO(\npixel\log(\npixel))$ operations using the FFT.
More details are given in section \ref{sec_Image_innerproducts_in_a_discrete_setting} below, and in \cite{Barnett2017,RSAB20}.

\subsection{Image discretization}
\label{sec_Image_discretization}

We denote by $\Omega_{1}$ and $\Omega_{\kmax}$ the ball of radius $1$ and $\kmax$, respectively, in either real- or frequency-space:
\begin{equation}
\Omega_{1} := \{ \vx \in \Real^{2} \mbox{~such~that~} \|\vx\| \leq 1 \} \text{,} \quad \Omega_{\kmax} := \{ \vk \in \Real^{2} \mbox{~such~that~} \|\vk\| \leq \kmax \} \text{,}
\end{equation}
and we will assume that all the images considered are supported in $\vx\in\Omega_{1}$.
Given that $A(\vx)$ is supported on $\Omega_{1}$, the representation $\fA(\vk)$ will have a bandlimit of $1$, implying that $\fA$ can be accurately reconstructed from its values sampled on a frequency-grid with spacing $\bigO(1)$ \cite{finufft}.

We also assume that any relevant signal within the images has a maximum effective spatial-frequency magnitude of $\kmax$; i.e., that the salient features of $\fA$ are concentrated in $\vk\in\Omega_{\kmax}$.
Consequently, we expect that the inversion
\begin{equation}
A(\vx) \approx \frac{1}{(2\pi)^{2}} \iint_{\Omega_\kmax} \fA(\vk)\euler^{+\imunit \vk \cdot \vx} \mathop{d\vk}
\end{equation}
will hold to high accuracy.
When these assumptions hold we can be sure that we won't lose much accuracy when applying these transformations in a discrete setting (as discussed below).

As discussed in \cite{RSAB20}, an $\npixel\times\npixel$ image spanning $[-1,+1]^{2}$ has a pixel-spacing of $\dx = 2/\npixel$, corresponding to a Nyquist spatial-frequency of $\kmax = \pi/\dx = (\pi/2)\npixel$.
We'll assume that the maximum effective spatial-frequency $\kmax$ is always taken to be on the order of the Nyquist spatial-frequency.
Thus, $\kmax = O(\npixel)$, and scalings in $\kmax$ and $\npixel$ will be equivalent.
Additionally, the fact that $x\leq 1$ implies that the bessel-coefficients $\bA(k,q)$ will be concentrated in the range $|q|\lesssim k$, meaning that the bessel-coefficients $\bA$ across all $k\in[0,\kmax]$ will be concentrated in $q\in[-\qmax/2,+\qmax/2-1]$ for $\qmax=\bigO(\kmax)=\bigO(\npixel)$.

With the notation above, we can consider an $\npixel\times\npixel$ image as a discrete set of pixel-averaged samples within $[-1,+1]^{2}$:
\begin{equation}
A_{n_{1},n_{2}} = \frac{1}{\dx^{2}}\int_{\vx_{1}=n_{1}\dx}^{(n_{1}+1)\dx} \int_{\vx_{2}=n_{2}\dx}^{(n_{2}+1)\dx} A(\vx) \mathop{d\vx}
\end{equation}
for indices $n_{1},n_{2}\in \{0,\ldots,\npixel-1\}$.
We approximate the fourier transform $\fA$ at any $\vk \in \Real^{2}$ via the simple summation:
\begin{equation}
\fA(\vk) = (\dx)^{2} \sum_{n_{2}=0}^{\npixel-1} \sum_{n_{1}=0}^{\npixel-1} A_{n_{1},n_{2}} \exp\left(-\imunit \vk \cdot \vx_{n_{1},n_{2}}\right)~ \text{,}
\label{eq_Ahattrap}
\end{equation}
where $\vx_{n_{1},n_{2}}$ is the appropriately-chosen pixel-center $\dx\left(n_{1}+\frac{1}{2},n_{2}+\frac{1}{2}\right)$.
Because we have assumed that the image is sufficiently well sampled (i.e., that $\fA$ contains little relevant frequency-content above the Nyquist-frequency $\kmax$), we expect the simple sum above to be accurate.

We will typically evaluate $\fourier{A}(\vk)$ for $\vk$ on a polar-grid, with $k$- and $\psi$-values corresponding to a quadrature-scheme, using the NUFFT to compute $\fA(k,\psi)$ at the associated quadrature-nodes (see \cite{finufft,RSAB20}).
As an example, for the application of image-alignment we use a Gauss-Jacobi quadrature for $k$ built with a weight-function corresponding to the radial-weighting of $kdk$.
This quadrature-scheme produces a set of $\rmax$ quadrature-nodes $k_{1},\ldots,k_{\rmax}$ and radial weights $w_{1},\ldots,w_{\rmax}$ which can be used to approximate the integral
\begin{equation}
\int_{0}^{\kmax} g(k)kdk \approx \sum_{r=1}^{\rmax} g(k_{r})w_{r} 
\end{equation}
to high accuracy for any function $g(k)$ which is smooth on the scale of $\bigO(1)$.
The $\qmax$ angular-nodes $\psi_{0},\ldots,\psi_{\qmax-1}$ will then be equispaced in the periodic interval $[0,2\pi)$, with a spacing of $\dpsi=2\pi/\qmax$, and $\psi_{q'}=q'\dpsi$.
Equispaced $\psi$-nodes allow for spectrally-accurate trapezoidal-quadrature in the $\psi$-direction, and we approximate the fourier-bessel-coefficients of each image-ring $\fA(k,\psi)$ as follows:
\begin{equation}
\bA(k,q) \approx \sum_{q^{\prime}=0}^{\qmax-1} \fA(k,\psi_{q^{\prime}}) \exp\left(-\imunit q \psi_{q^{\prime}}\right) \dpsi \text{,}
\end{equation}
with the index $q$ considered periodically in the interval $[-\qmax/2-1,\ldots,+\qmax/2]$ (so that, e.g., the $q$-value of $\qmax-1$ corresponds to the $q$-value of $-1$).
%When dealing solely with images in 2-dimensions, we typically use a Gauss-Jacobi quadrature built with Jacobi-polynomials of type $\alpha=0$, $\beta=1$, which is optimized for the radial-weighting of $kdk$.
%When dealing with images that are projections of a 3d-volume, we often use a Gauss-Jacobi quadrature built with Jacobi-polynomials of type $\alpha=0$, $\beta=2$.
%While this quadrature-scheme is optimized for the radial-weighting of $k^{2}dk$ in 3-dimensions, we can easily solve a linear system to construct specialized weights which accommodate the radial-weighting of $kdk$ required for integration in 2-dimensions.

\begin{figure}
\centering
\includegraphics[width=7in]{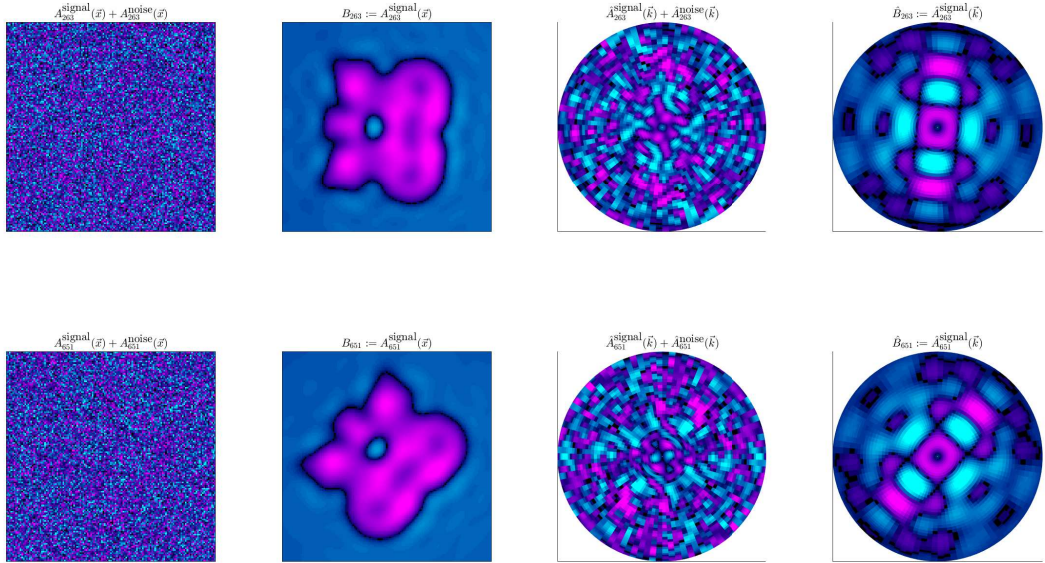}
\caption{\label{pm_fig_image_pair_FIGA}
Images of the TRPV1 molecule taken from \cite{Liao2013}.
From left to right we illustrate $A(\vx)$, $A^{\signal}(\vx)$, $\fA(k,\psi)$ and $\fA^{\signal}(k,\psi)$ for one image top, with analogous representations of another image on the bottom.
Real-space images are shown in a box of $137\times 137$ pixels, corresponding to a $165$\r{A} side-length square mask.
Fourier-space images (real-part only) are shown in a polar-grid with $\rmax=49$ radial quadrature-nodes with $\kmax=48$ and $\dpsi = 2\pi/98$.
In each case the colormap is centered on the mean (black) and ranges across $\pm 2.5$ standard-deviations (magenta to cyan).
}
\end{figure}

An example of some discretized images is given in Fig \ref{pm_fig_image_pair_FIGA}.
This figure illustrates images of the TRPV1 molecule taken from the Electron Microscopy Public Image Archive dataset {\tt EMPIAR-10005} \cite{Liao2013}.
The top and bottom rows correspond to two different images, indexed by their order in the dataset.
The original images $A(\vx)$ (far left) are quite noisy, and we estimate their `true' signals $A^{\signal}$ (shown adjacent) as follows.
First, we project the 3-dimensional electron-density-function of the TRPV1-reference-molecule (taken from the Electron Microscopy Data Bank (EMDB) structure `{\tt EMD-5778}') onto the estimated viewing-angles of the images (in $SO3)$.
These projections form `templates', which we then convolve with the fourier-transform of the contrast-transfer-function (CTF) associated with the imaging process.
The two images have been chosen so that their estimated viewing-angles in $SO3$ are roughly the same, up to an in-plane rotation of $\sim 37^{\circ}$.
Thus, the two images have approximately the same template (modulo in-plane rotation), which we refer to as $S(\vx)$.
These two images also have approximately the same CTF.
Thus, the signal for both images can be estimated using the same `CTF-corrected template' $CTF\odot\fS$.
More specifically, we approximate the image signals $\fB:=\fA^{\signal}$ as $\rotation_{\gamma}\left[CTF\odot\fS\right](\vx)$ for some appropriately chosen image-specific in-plane angle $\gamma$.
Representations of $\fA(k,\psi)$ and $\fA^{\signal}(k,\psi)$ are shown on the right side of Fig \ref{pm_fig_image_pair_FIGA}.
These are displayed on a polar-grid out to $\kmax=48$, which is a typical frequency associated with low-resolution analysis and processing in this setting.

\begin{figure}
\centering
\includegraphics[width=6in]{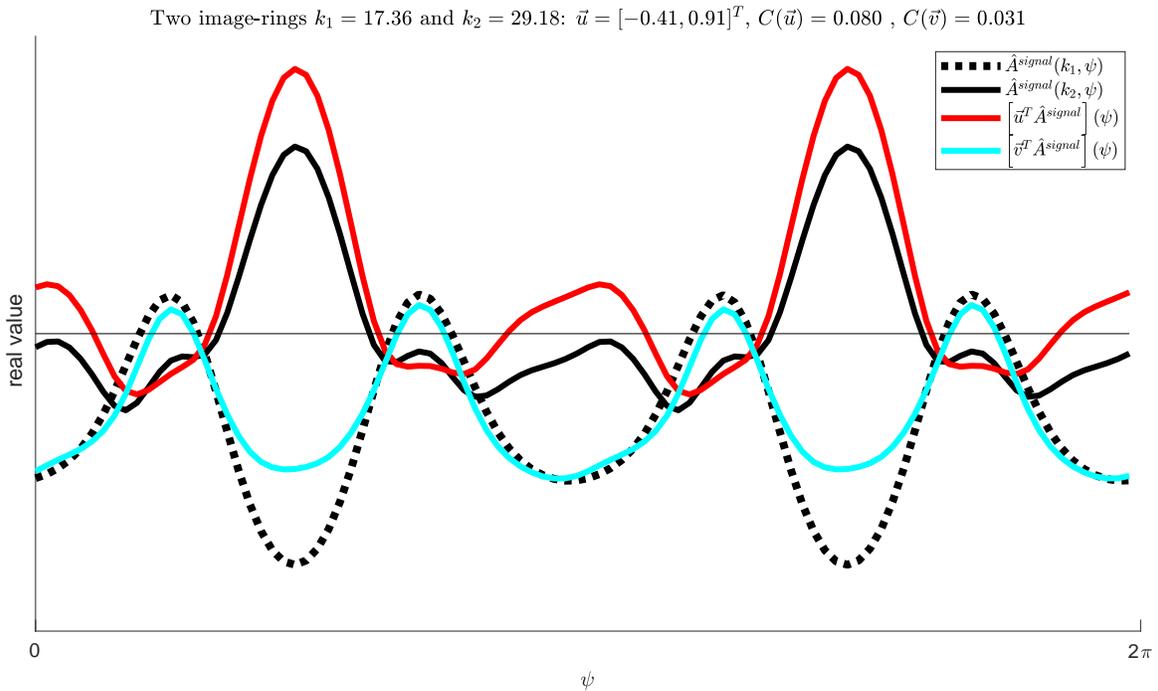}
\caption{\label{pm_fig_image_pair_S2_FIGB}
Here we show two different image-rings $\fA^{\signal}(k_{1},\psi)$ and $\fA^{\signal}(k_{2},\psi)$ in dashed- and solid-black lines, respectively.
The principal-vector $\vu$ is selected to maximize the simple objective-function $\Ccost$ described in the text.
The orthogonal vector $\vv$ minimizes $\Ccost$.
The dominant principal-image-ring $\transpose{\vu}\fA$ is shown in red, and $\transpose{\vv}\fA$ is shown in cyan.
}
\end{figure}

\begin{figure}
\centering
\includegraphics[width=6in]{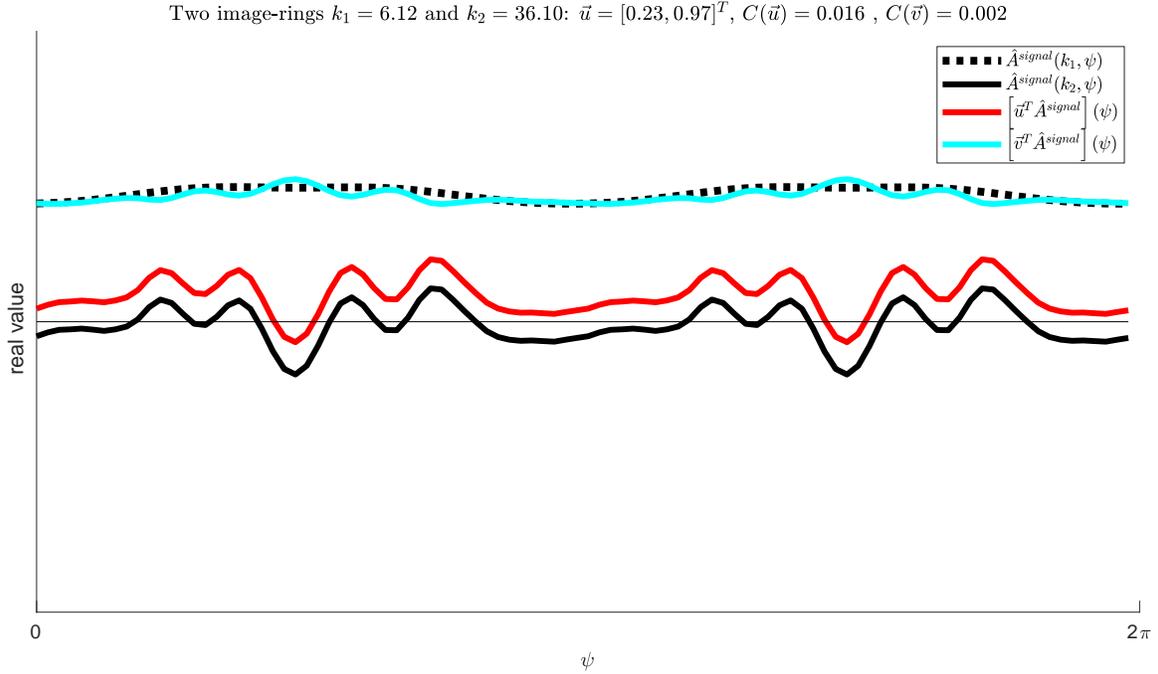}
\caption{\label{pm_fig_image_pair_S2_FIGC}
This figure has the same format as Fig \ref{pm_fig_image_pair_S2_FIGB}, for a different pair of image-rings.
}
\end{figure}

\begin{figure}
\centering
\includegraphics[width=3.5in]{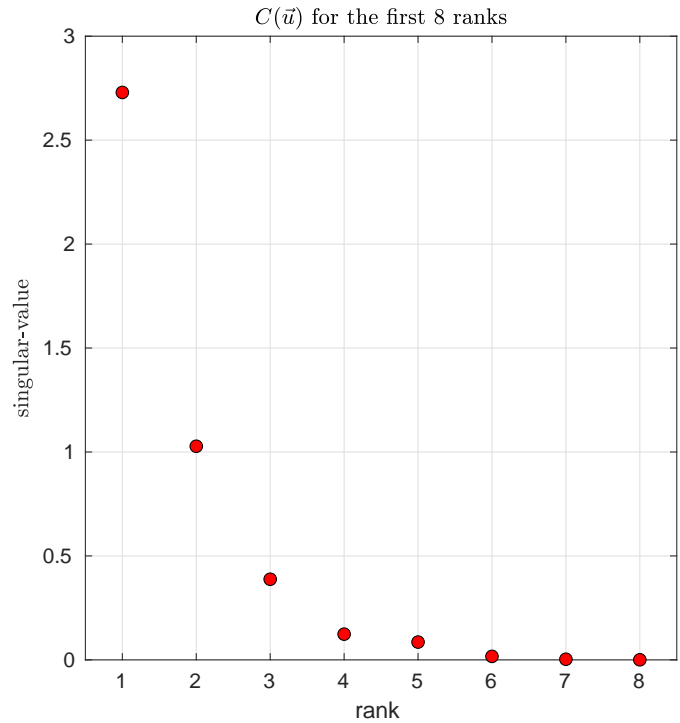}
\caption{\label{pm_fig_image_pair_Sall_FIGD}
Here we show the eigenvalues of $\Ckern$ for the objective-function $\Ccost(\vu)$ described in the text.
Note that, while the full rank of $\Ckern$ is $\rmax=49$, the spectrum of $\Ckern$ has decayed considerably by $\nrank=6-7$.
}
\end{figure}

\begin{figure}
\centering
\includegraphics[width=7in]{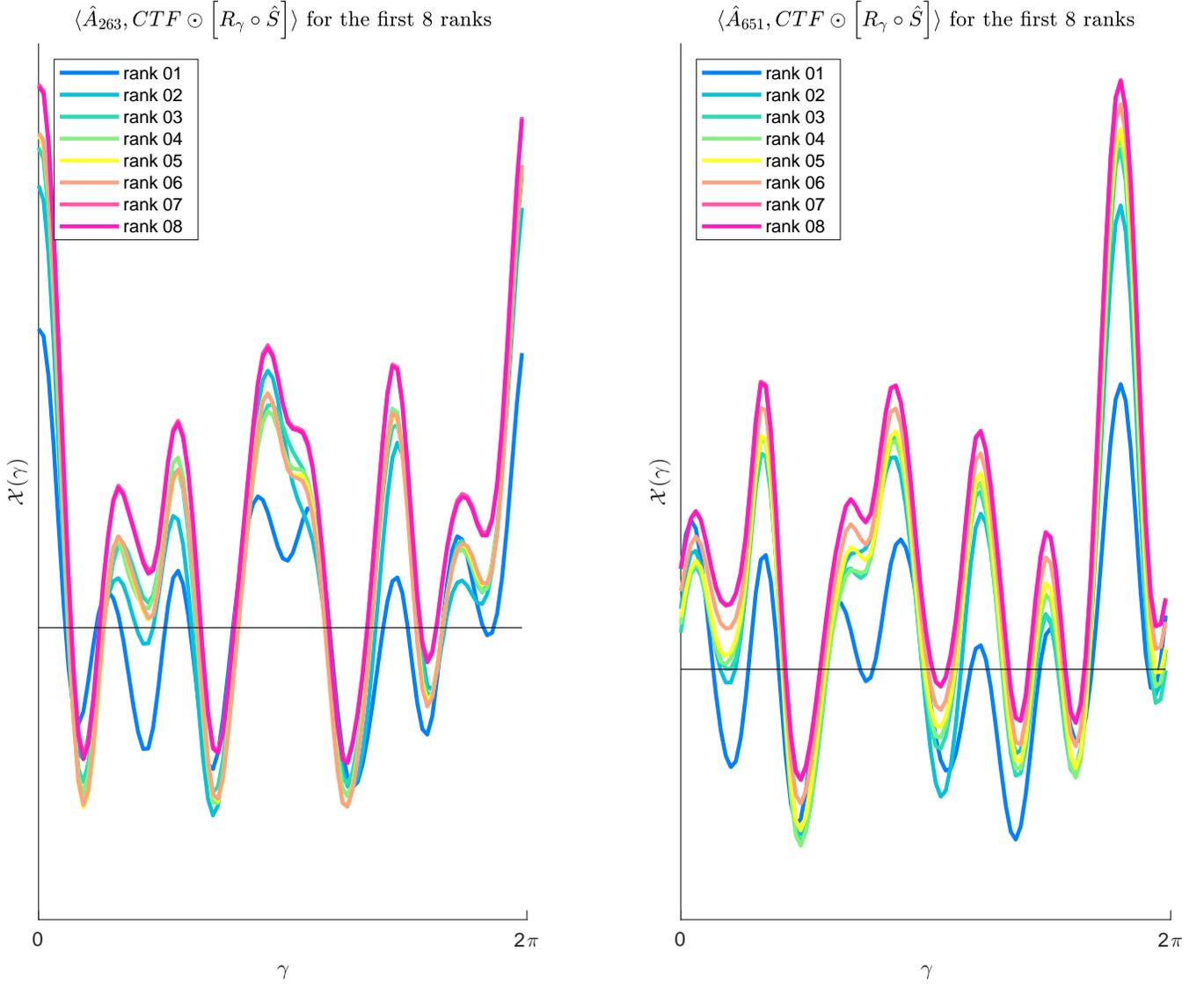}
\caption{\label{pm_fig_image_pair_X_wSM_FIGE}
Here we show the inner-product landscape $\cX(\gamma ; \fA,CTF\odot\fS )$ as a function of $\gamma$, for both of the images shown in Fig \ref{pm_fig_image_pair_FIGA}.
In this expression the `CTF' refers to the contrast-transfer-function.
The inner-products are calculated over a range of $\nrank\in\left\{1,\ldots,8\right\}$.
Note that, even when $\nrank$ is low and the function $\cX(\gamma)$ is inaccurate, the optimal alignment-angle is usually still accurate.
}
\end{figure}

\begin{figure}
\centering
\includegraphics[width=3.5in]{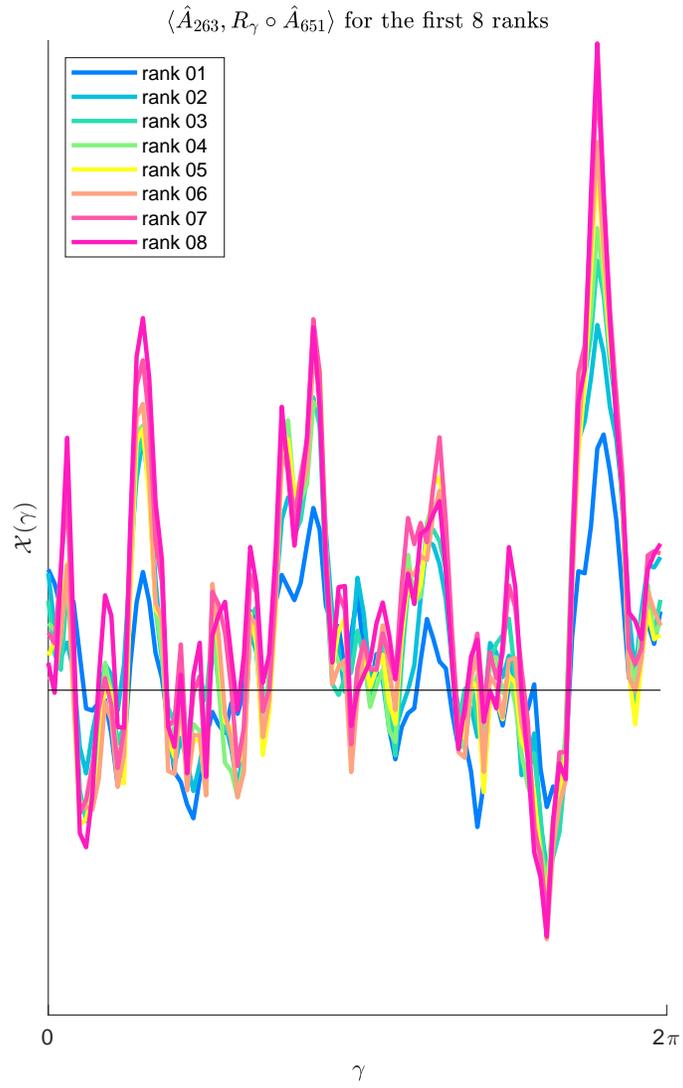}
\caption{\label{pm_fig_image_pair_X_wSM_FIGF}
This figure has the same format as Fig \ref{pm_fig_image_pair_X_wSM_FIGF}, except that it compares the two experimental images to one another.
}
\end{figure}

\begin{figure}
\centering
\includegraphics[width=7in]{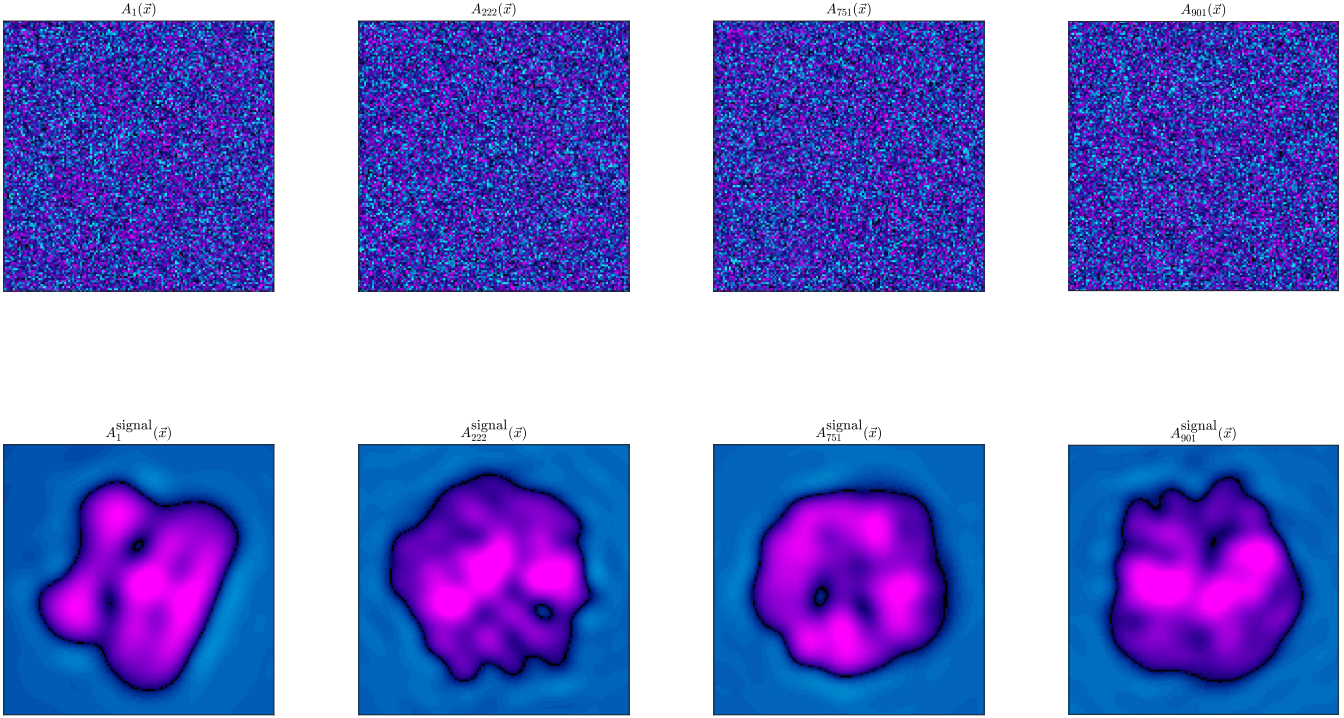}
\caption{\label{pm_fig_M_and_S_FIGG}
Here we show several images (top) from the Electron Microscopy Public Image Archive dataset {\tt EMPIAR-10005}.
The associated CTF-corrected-templates are shown below.
The format for these images is the same as Fig \ref{pm_fig_image_pair_FIGA}.
}
\end{figure}

\begin{figure}
\centering
\includegraphics[width=7in]{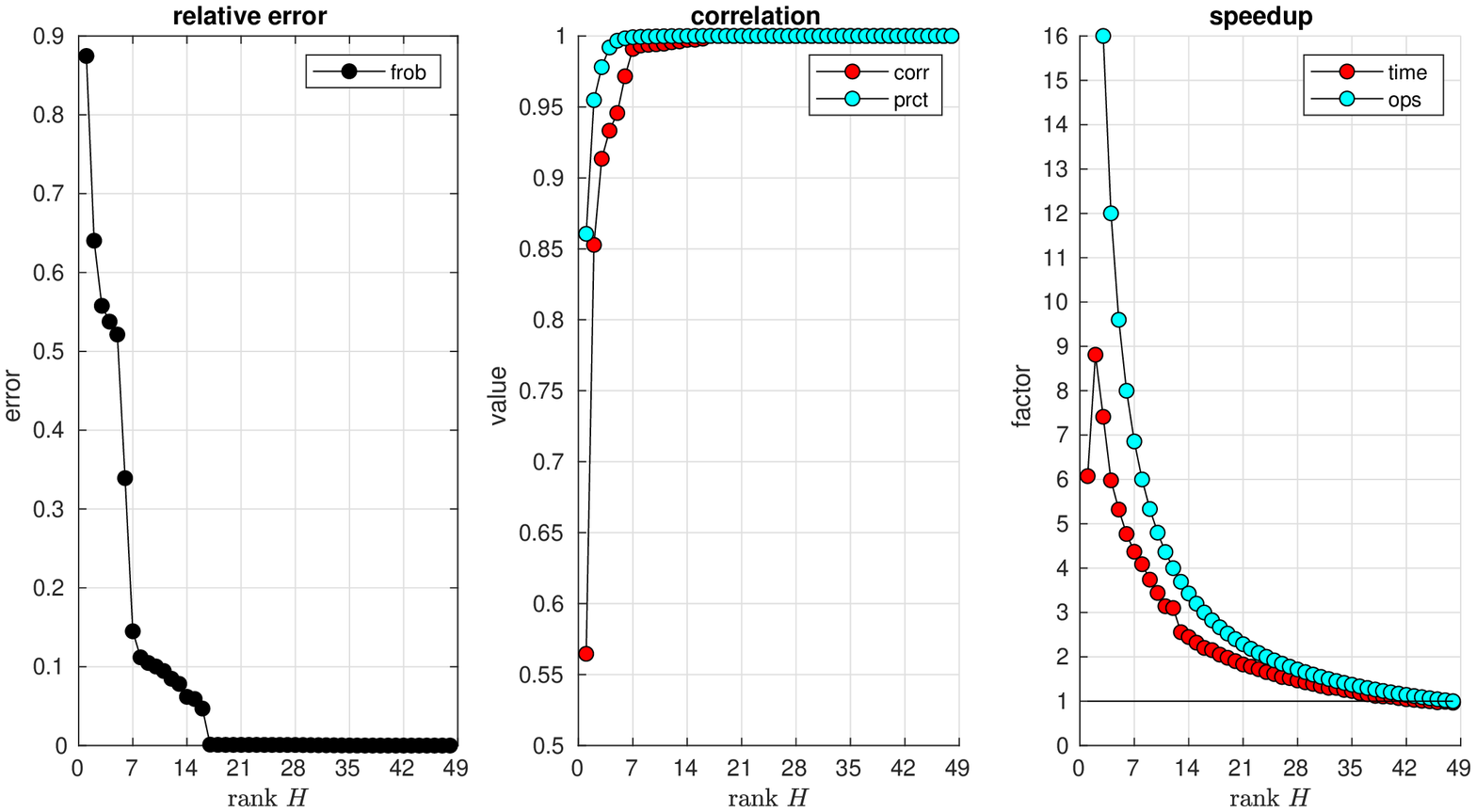}
\caption{\label{pm_fig_X_2d_cluster_d0_FIGK}
Here we show the results after applying our radial-SVD to align the first $1024$ images in the {\tt EMPIAR-10005} dataset to a set of $993$ targets (see main text).
We use the strategy of section \ref{sec_Image_innerproducts_in_a_discrete_setting} to calculate $\cX(\gamma ; \fA,\fB )$ for each image-target pair $(A,B)$.
We approximate this inner-product array by using the radial-SVD, retaining $\nrank$ principal-vectors of $\Ckern$ as described in the main text.
For each $\nrank$, we calculate $\cX^{\estimated}(\gamma ; \fA,\fB ; \nrank)$, and compare the results to the full calculation (involving all $\rmax$ radial quadrature-nodes).
The relative error (using the frobenius-norm) between these two arrays is shown in black on the left.
We emphasize that the frobenius-norm does not tell the whole story; our objective-function is designed to emphasize alignment accuracy, not the magnitude of the inner-product.
This is reflected in the middle subplot, which shows two other measures of alignment accuracy (see main text).
In the limit as the number of image-target pairs becomes very large, we expect a speedup of roughly $\rmax/\nrank$, as shown in cyan in the right subplot.
The actual speedup in total runtime (which includes the necessary precomputations), for an implementation on a dell laptop with an i7 processor, is shown in red.
}
\end{figure}

\subsection{Image inner-products in a discrete setting}
\label{sec_Image_innerproducts_in_a_discrete_setting}

As alluded to in section \ref{sec_Image_innerproducts_in_a_continuous_setting}, we can use the discrete fourier-bessel-coefficients $\bA(k,q)$ and $\bB(k,q)$ associated with images $\fA$ and $\fB$ to efficiently calculate the inner-products $\cX$ across a range of rotation-angles $\gamma_{q'}=2\pi q'/\qmax$.
This calculation can be summarized as:
\begin{equation}
\label{eq_image_innerproduct_0}
\hspace{-15mm}
\cX(\gamma_{q^{\prime}} ; \fA,\fB) = 2\pi \sum_{q=0}^{\qmax-1} \exp\left(-\imunit q \gamma_{q^{\prime}} \right) \sum_{r=1}^{\rmax} w_{r} \bA(k_{r},q)^{\dagger} \bB(k_{r},q) \text{.}
\end{equation}
This calculation can be achieved in two steps:
\begin{flalign}
\label{eq_image_innerproduct}
& \textbf{[Step 1]} \quad \fourier{\cX}(q ; \fA,\fB) = 2\pi \sum_{r=1}^{\rmax} w_{r} \bA(k_{r},q)^{\dagger}\bB(k_{r},q) \text{,} && \\\nonumber
& \textbf{[Step 2]} \quad \cX(\gamma_{q^{\prime}} ; \fA,\fB) = \sum_{q=0}^{\qmax-1} \exp\left(-\imunit q \gamma_{q^{\prime}}\right) \fourier{\cX}(q ; \fA,\fB) \text{.} &&
\end{flalign}
The first step applies radial-quadrature to combine information from different image-rings, computing the fourier series coefficients $\fourier{\cX}(q)$ in $\bigO(\rmax \qmax)$ operations.
The second step can be evaluated using a 1-dimensional fast-fourier-transform of size $\qmax$, requiring $\bigO(\qmax \log(\qmax))$ operations.
If a finer resolution in $\gamma$ is required then the array $\fourier{\cX}$ can be zero-padded prior to the FFT.
If a specific non-uniform set of $\gamma$-values are requested, then then non-uniform fast-fourier-transform (NUFFT) can be applied instead.

In most cases the operation-count is dominated by the first step -- i.e., the calculation of the $\qmax$ values $\fourier{\cX}(q)$.
This operation-count scales with $\rmax$, which is the number of image-rings (i.e., distinct $k$-values) in our polar quadrature-grid.
The power-spectral-density of an arbitrary image-array $A_{n_{1},n_{2}}$ is typically constant out to $\vk\leq\kmax$, implying that an accurate calculation of $\fourier{\cX}$ will generally require all $\rmax$ of these radial quadrature-nodes.
However, as we will emphasize below, images that arise within the context of cryo-EM are far from arbitrary, and have frequency-content which can usually be exploited to accelerate the calculation of $\fourier{\cX}$.

In the simplest situation, one might imagine an image with all the relevant frequency content restricted to a single value of $k_{r}$, corresponding to a single image-ring:
\begin{equation}
\fA(\vk) = \begin{cases}
    \fA(k_{r},\psi), & \text{if } k=k_{r} \\
    0,              & \text{otherwise.}
\end{cases}
\end{equation}
In this idealized case one could consider only the single summand for $k_{r}$ in \eqref{eq_image_innerproduct}-1, thus calculating $\fourier{\cX}$ with only $\bigO(\qmax)$ operations, saving a factor of $\rmax$.

While real images are never quite this simple, they often allow for a similar reduction in operation-count.
The strategy we will discuss later on in section \ref{sec_Radial_svd} involves restructuring the image by taking linear-combinations of different image-rings.
In order to explain which linear-combinations we choose (and why), we first review the standard noise-model for the images.
As described later on in section \ref{sec_Image_objective_single}, we will use this noise-model to motivate a simple quadratic objective-function to inform a rank-reducing transformation which compresses the frequency-content of the images.
After this compression, fewer than $\rmax$ terms will be needed to accurately compute $\fourier{\cX}$.
While this strategy may not amount to much of a computational savings when performing a single image-image alignment, it will be quite useful when aligning the same image multiple times to  multiple targets, which is often done over the course of cryo-EM molecular-reconstruction.

\subsection{Image noise model}
\label{sec_Image_noise_model}

Because of our focus on cryo-EM molecular-reconstruction, we will assume that the images considered can be modeled as the sum of a `signal' plus a `noise'.
We assume that the signal corresponds to a 2-dimensional projection of a (smooth) molecular electron-density-function, while the noise corresponds to detector noise \cite{Sigworth1998,Scheres2009,Sigworth2010,Lyumkis2013}.
We do not consider more complicated sources of noise, such as structural-noise associated with image-preparation (see \cite{Baxter2009}).

Based on these simple assumptions, we'll model the noise in real-space as independent and identically-distributed (iid), with a variance of $\sigma^{2}$ on a unit-scale in 2-dimensional real-space.
In this section and the next, we review this standard noise-model using the notation we've introduced above.

To begin with, we'll model the discrete image values
\begin{equation}
A_{n_{1},n_{2}} = A(\vx_{n_{1},n_{2}}) = A^{\signal}_{n_{1},n_{2}} + A^{\noise}_{n_{1},n_{2}} \text{,} 
\end{equation}
where the signal and noise are represented by the arrays $A^{\signal}_{n_{1},n_{2}}$ and $A^{\noise}_{n_{1},n_{2}}$, respectively.
Because each entry in $A^{\noise}_{n_{1},n_{2}}$ corresponds to an average over the area-element $\dx^{2}$ associated with a single pixel, we expect the variance of each $A^{\noise}_{n_{1},n_{2}}$ to scale inversely with $\dx^{2}$.
That is to say, we'll assume that each element of the noise-array $A^{\noise}_{n_{1},n_{2}}$ is drawn from the standard normal distribution $\cN\left(0,\frac{\sigma^{2}}{\dx^{2}}\right)$, with a mean of $0$ and a variance of $\sigma^{2}/\dx^{2}$.
With these assumptions, downsampling the image by averaging neighboring pixels will correspond to an increase in the pixel-size and a simultaneous reduction in the variance of the noise associated with each (now larger) pixel.

Given the assumptions above, we expect that $\fA$ will be modeled by:
\begin{equation}
\fA(\vk) = \fA(\vk) = \fA^{\signal}(\vk) + \fA^{\noise}(\vk) \text{,} 
\end{equation}
where $\fA^{\noise}$ is now complex and iid.
Because the noise-term $A^{\noise}$ is real, the noise-term $\fA^{\noise}$ will be complex, with the conjugacy-constraint that $\fA^{\noise}(+\vk) = \fA^{\noise}(-\vk)^{\dagger}$.

If we were to sample $\vk$ on the uniform $\npixel\times\npixel$ cartesian-grid associated with the standard 2-dimensional fast-fourier-transform, then $\vk_{n_{1},n_{2}}$ will correspond to the frequency $\dk(n_{1},n_{2})$, with frequency-spacing $\dk=\pi$ and indices $(n_{1},n_{2})$ considered periodically in the range $[-\npixel/2,\ldots,+\npixel/2-1]$.
In this case the transformation between $A_{n_{1},n_{2}}$ and $\fA_{n_{1},n_{2}}$ will be unitary.
When $(n_{1},n_{2})$ corresponds to an index-pair that is periodically reflected onto itself (i.e., the four index-pairs $(0,0)$, $(\npixel/2,0)$, $(0,\npixel/2)$ and $(\npixel/2,\npixel/2)$) then $\fA^{\noise}_{n_{1},n_{2}}$ will be real and drawn from $\cN(0,\frac{\sigma^{2}}{\dx^{2}})$.
When $(n_{1},n_{2})$ corresponds to an index-pair that is {\em not} periodically reflected onto itself, the real- and imaginary-components of $\fA^{\noise}(\vk_{n_{1},n_{2}})$ will each be drawn iid from $\cN(0,\frac{\sigma^{2}}{2\dx^{2}})$ subject to the conjugacy constraint above.% (note the additional factor of 2 in the variance).

Motivated by this observation, we define $\fsigma^{2}=\frac{\pi^{2}\sigma^{2}}{\dx^{2}}$ as the variance on a unit-scale in 2-dimensional frequency-space.
We expect that the noise-term $\fA^{\noise}(\vk)$ integrated over any area-element $\dk^{2}$ in frequency-space will have a variance of $\frac{\fsigma^{2}}{\dk^{2}}$.
Recalling our polar quadrature described above, we typically record $\vk$ along $\rmax$ radial quadrature-nodes $k_{1},\ldots,k_{\rmax}$ and $\qmax$ angular quadrature-nodes $\psi_{0},\ldots,\psi_{\qmax-1}$, with radial- and angular-weights $w_{r}$ and $\dpsi$, respectively.
In this quadrature-scheme each quadrature-node $\vk_{rq}=(k_{r},\psi_{q})$ is associated with the area-element $w_{r}\dpsi$, which approximates the `$kdk$' integration weight.
Consequently, we expect that the noise-term $\fA^{\noise}(\vk_{rq})$ evaluated at $(k_{r},\psi_{q})$ on our polar quadrature-grid will have a variance of $\frac{\fsigma^{2}}{w_{r}\dpsi}$.

\subsection{Image similarity: likelihood}
\label{sec_Image_similarity}

Given one `noisy' image $\fA=\fA^{\signal}+\fA^{\noise}$ and another `noiseless' image $\fB=\fB^{\signal}$, we can derive the standard likelihood $P(\fA \given[\big] \fA^{\signal}=\fB)$ of observing any value of $\fA$, given that the signals are the same.

We first note that, for a particular quadrature-node $\vk_{rq} = (k_{r},\psi_{q})$ on our polar-grid, the image-value $\fA(\vk_{rq})$ will be drawn from the gaussian-distribution:
\begin{equation}
\fA(\vk_{rq}) \sim \cN(\fA^{\signal}(\vk_{rq}),\frac{\fsigma^{2}}{w_{r}\dpsi}) \text{.}
\end{equation}
Generally speaking, the variance of this distribution will depend on the choice of the image-ring (i.e., the index $r$).
To match the variance of these distributions (across image-rings), we consider the rescaled image-values $\fA(\vk_{rq})\eta_{r}$ and $\fB(\vk_{rq})\eta_{r}$, with $\eta_{r}=\sqrt{w_{r}\dpsi}$
The rescaled image-values will be drawn from the gaussian-distribution:
\begin{equation}
\fA(\vk_{rq})\eta_{r} \sim \cN(\fA^{\signal}(\vk_{rq})\eta_{r},\fsigma^{2}) \text{,}
\end{equation}
with a variance that no longer depends on the choice of image-ring.

The likelihood of any observation of $\fA(\vk)$, given the hypothesis that the signal $\fA^{\signal}(\vk)$ equals $\fB(\vk)$ can be derived from the relationship:
\begin{equation}
P(\fA^{\signal}=\fB \given[\big] \fA) \cdot P(\fA) = P(\fA \given[\big] \fA^{\signal}=\fB) \cdot P(\fA^{\signal}=\fB) \text{,}
\end{equation}
which can be rearranged into:
\begin{equation}
P(\fA^{\signal}=\fB \given[\big] \fA) = \frac{P(\fA \given[\big] \fA^{\signal}=\fB)\cdot P(\fA^{\signal}=\fB)}{P(\fA)} \text{.}
\end{equation}

Assuming a uniform prior $P(\fA^{\signal}=\fB)$, the posterior-probability $P(\fA^{\signal}=\fB|\fA)$ is proportional to the likelihood $P(\fA \given[\big] \fA^{\signal}=\fB)$.
For a single quadrature-node $(k_{r},\psi_{q})$ this likelihood is given by:
\begin{eqnarray}
P(\fA(\vk_{rq})\eta_{r} \given[\big] \fA^{\signal}(\vk_{rq})\eta_{r}=\fB(\vk_{rq})\eta_{r} ; r , q \text{\ fixed}) \quad = \\
\qquad \frac{1}{\sqrt{2\pi}\fsigma} \exp\left(- \frac{ | \fA(\vk_{rq}) - \fB(\vk_{rq}) |^{2} }{ 2\fsigma^{2} } w_{r}\dpsi \right) \text{.}
\end{eqnarray}
Because the values of $\fA^{\noise}(\vk_{rq})$ are independent from one another, the likelihood of observing the image-ring $\fA(\vk_{rq})$ given the hypothesis that one image-ring of $\fA^{\signal}(\vk_{rq})\eta_{q}$ is equal to the corresponding image-ring $\fB(\vk_{rq})\eta_{q}$ (taken across $q$) is:
\begin{eqnarray}
P(\fA(\vk_{rq})\eta_{r} \given[\big] \fA^{\signal}(\vk_{rq})\eta_{r}=\fB(\vk_{rq})\eta_{r} ; r \text{\ fixed}) \quad = \\
\qquad \Pi_{q=0}^{\qmax-1} P(\fA(\vk_{rq})\eta_{r} \given[\big] \fA^{\signal}(\vk_{rq})\eta_{r}=\fB(\vk_{rq})\eta_{r} ; r , q \text{\ fixed}) \text{.}
\end{eqnarray}
Taking the logarithm of this expression, we see that the log-likelihood is (up to an additive constant of $\qmax\log(\sqrt{2\pi}\fsigma)$, which we ignore):
\begin{eqnarray}
\log P(\fA(\vk_{rq})\eta_{r} \given[\big] \fA^{\signal}(\vk_{rq})\eta_{r}=\fB(\vk_{rq})\eta_{r} ; r \text{\ fixed}) \quad = \\
\qquad - \frac{1}{2\fsigma^{2}} \sum_{q=0}^{\qmax-1} |\fA(\vk_{rq}) - \fB(\vk_{rq})|^{2} w_{r}\dpsi \text{,}
\end{eqnarray}
which converges to
\begin{eqnarray}
\log P(\fA(\vk_{rq})\eta_{r} \given[\big] \fA^{\signal}(\vk_{rq})\eta_{r}=\fB(\vk_{rq})\eta_{r} ; r \text{\ fixed}) \quad = \\
\qquad - \frac{1}{2\fsigma^{2}} w_{r}\cdot \int_{0}^{2\pi} |\fA(k_{r},\psi) - \fB(k_{r},\psi)|^{2} d\psi \text{\ }
\label{eq_L_B_A_ring}
\end{eqnarray}
as the number of quadrature-nodes $\qmax\rightarrow\infty$.
Similarly, the likelihood of observing the array $\fA(\vk_{rq})$, given the hypothesis that the array $\fA^{\signal}(\vk_{rq})$ is equal to the array $\fB(\vk_{rq})\eta_{r}$ (taken across $r,q$) is:
\begin{eqnarray}
P(\fA(\vk_{rq})\eta_{r} \given[\big] \fA^{\signal}(\vk_{rq})\eta_{r}=\fB(\vk_{rq})\eta_{r}) \quad = \\
\qquad \Pi_{r=1}^{\rmax} P(\fA(\vk_{rq})\eta_{r} \given[\big] \fA^{\signal}(\vk_{rq})\eta_{r}=\fB(\vk_{rq})\eta_{r} ; r \text{\ fixed}) \text{,}
\end{eqnarray}
corresponding to a log-likelihood of:
\begin{eqnarray}
\log P(\fA(\vk_{rq})\eta_{r} \given[\big] \fA^{\signal}(\vk_{rq})\eta_{r}=\fB(\vk_{rq})\eta_{r}) \quad = \\
\qquad - \frac{1}{2\fsigma^{2}} \sum_{r=1}^{\rmax}\sum_{q=0}^{\qmax-1} |\fA(\vk_{rq}) - \fB(\vk_{rq})|^{2} w_{r}\dpsi \text{,}
\end{eqnarray}
which converges to
\begin{eqnarray}
\log P(\fA(\vk_{rq})\eta_{r} \given[\big] \fA^{\signal}(\vk_{rq})\eta_{r}=\fB(\vk_{rq})\eta_{r}) \quad = \\
\qquad - \frac{1}{2\fsigma^{2}} \iint_{\Omega_{\kmax}} |\fA(\vk) - \fB(\vk)|^{2} d\vk \text{\ }
\label{eq_L_B_A_image}
\end{eqnarray}
as the number of quadrature-nodes $\rmax\rightarrow\infty$.
This final expression is exactly the standard result one would expect given iid-noise with a fixed variance on a unit-scale in frequency-space \cite{Sigworth1998,Scheres2009,Sigworth2010,Lyumkis2013}.

When dealing with two noisy images $\fA=\fA^{\signal}+\fA^{\noise}$ and $\fB=\fB^{\signal}+\fB^{\noise}$, the same general argument applies, so long as $\fA^{\noise}$ is independent from $\fB^{\noise}$.
The only difference is that now the value $\fB(\vk_{qm})\eta_{r}$ for each quadrature-node will not be fixed, but will be drawn from $\cN(\fB(\vk_{qm})\eta_{r},\fsigma^{2})$.
Consequently, the log-likelihood of the observations $\fA$ and $\fB$, given the hypothesis that the array $\fA^{\signal}$ is equal to $\fB^{\signal}$ is:
\begin{eqnarray}
\log P(\fA , \fB \given[\big] \fA^{\signal}=\fB^{\signal}) = - \frac{1}{4\fsigma^{2}} \iint_{\Omega_{\kmax}} | \fA(\vk) - \fB(\vk) |^{2} d\vk \text{.}
\end{eqnarray}
%Note that the variance has increased by a factor of $2$ -- i.e., $\cN(0,\sigma_{1}^{2})\star\cN(0,\sigma_{2}^{2})=\cN(0,\sigma_{1}^{2}+\sigma_{2}^{2})$.
Note that, aside from a constant factor, this is identical to \eqref{eq_L_B_A_image}.

\subsection{Image alignment: maximum-likelihood estimate}
\label{sec_Image_alignment}

We now review the standard approach for aligning two images $\fA=\fA^{\signal}+\fA^{\noise}$ and $\fB=\fB^{\signal}+\fB^{\noise}$ to one another.
If we assume that $\fA^{\signal}$ is equal to $\rotation_{\gamma}\fB^{\signal}$ for some unknown angle $\gamma$, then the maximum-likelihood estimate for $\gamma$ is given by:
\begin{eqnarray}
\gamma_{\optimal} & = & \argmax_{\gamma} \log P(\fA^{\signal}=\rotation_{\gamma}\fB^{\signal} \given[\big] \fA , \rotation_{\gamma}\fB) \\
               & = & \argmax_{\gamma} - \iint_{\Omega_{\kmax}} | \fA(\vk) - \rotation_{\gamma}\fB(\vk) |^{2} d\vk \text{.}
\label{eq_Image_alignment_maximum_likelihood}
\end{eqnarray}

Note that, if we were to ignore the $\gamma$-independent terms in \eqref{eq_Image_alignment_maximum_likelihood}, we would obtain:
\begin{eqnarray}
\gamma_{\optimal} & = & \argmax_{\gamma} \langle \fA , \rotation_{\gamma}\fB \rangle \\
               & = & \argmax_{\gamma} \cX(\gamma) \text{,}
\label{eq_Image_alignment_innerproduct}
\end{eqnarray}
which is the reason why the inner-product $\cX(\gamma)$ is often used as a measurement of similarity in image-alignment \cite{Sigworth1998,Sigworth2010}.

\section{Radial-SVD}
\label{sec_Radial_svd}

In this section we introduce the radial-SVD, focusing on 2-dimensional images for ease of presentation.
We first present an idealized case-study which motivates the radial-SVD using an objective-function for a single image-target pair.
Then we generalize the objective-function to multiple image-target pairs and present an example.

\subsection{Image case-study: a pair of image-rings}
\label{sec_Image_case_study}

Given the standard noise-model and likelihood formula above, we now present a simple case-study that motivates this paper.

Imagine that we are given two different images $\fA=\fA^{\signal}+\fA^{\noise}$ and $\fB=\fB^{\signal}$ which are both observations of the {\em same} signal $\fA^{\signal}=\fB^{\signal}$, but with only the first image $\fA$ affected by noise (the same argument will hold if $\fB$ is noisy as well, provided that $\fA^{\noise}$ and $\fB^{\noise}$ are independent).
Furthermore, let's assume for simplicity that there are only two distinct $k$-values $k_{1}$ and $k_{2}$ in our polar quadrature-grid (i.e., $\rmax=2$) and that the quadrature is exact with unit weights (so $w_{1}=w_{2}=1$).
Thus, the two images are both restricted to only two image-rings corresponding to two distinct $k$-values $k_{1}$ and $k_{2}$.
The image $\fA$ comprises the two image-rings $\fA(k_{1},\psi)$ and $\fA(k_{2},\psi)$, each a function of $\psi$; similar notation holds for $\fB$.

We now consider the problem of aligning these images to one another.
That is, we assume that $\fA^{\signal}$ is equal to $\rotation_{\gamma}\fB^{\signal}$ for some unknown angle $\gamma$.
As reviewed above, the maximum-likelihood estimate for $\gamma$ is given by:
\begin{eqnarray}
\gamma_{\optimal} & = & \argmax_{\gamma} \log P(\fA^{\signal}=\rotation_{\gamma}\fB^{\signal} \given[\big] \fA , \rotation_{\gamma}\fB) \\
               & = & \argmax_{\gamma} - \iint_{\Omega_{\kmax}} | \fA(\vk) - \rotation_{\gamma}\fB(\vk) |^{2} d\vk \\
               & = & \argmax_{\gamma} - \int_{0}^{2\pi} | \fA(k_{1},\psi) - \fB(k_{1},\psi-\gamma) |^{2} + | \fA(k_{2},\psi) - \fB(k_{2},\psi-\gamma) |^{2} d\psi \text{,}
\label{eq_example_maximum_likelihood}
\end{eqnarray}
which, as shown in the last line, involves the sum of the log-likelihood associated with each image-ring.

If we have the computational resources available, then we can certainly calculate the log-likelihood for each image-ring, sum the results, and obtain the full log-likelihood.
However, if we do not have the resources available to perform the full calculation, or if we simply want to approximate the full log-likelihood, we might restrict our calculation to only one of the $k_{r}$.
That is, we might consider the approximate maximum-likelihood estimate:
\begin{eqnarray}
\gamma_{\optimal}^{\estimated} & = & \argmax_{\gamma} \log P(\fA^{\signal}=\rotation_{\gamma}\fB^{\signal} \given[\big] \fA , \rotation_{\gamma}\fB ; r \text{\ fixed}) \\
              & = & \argmax_{\gamma} - \int_{0}^{2\pi} | \fA(k_{r},\psi) - \fB(k_{r},\psi-\gamma) |^{2} d\psi \text{,}
\label{eq_example_maximum_likelihood_single_ring}
\end{eqnarray}
for either $r=1$ or $r=2$.
As described in section \ref{sec_Image_innerproducts_in_a_discrete_setting}, the operation-count associated with calculating \eqref{eq_example_maximum_likelihood_single_ring} will be lower than that for calculating the full likelihood in \eqref{eq_example_maximum_likelihood}, simply because \eqref{eq_example_maximum_likelihood_single_ring} involves only a single radial quadrature-node, whereas \eqref{eq_example_maximum_likelihood} involves two.

If we insist on such an approximation, the natural question is: which of the two image-rings is a better choice to use in \eqref{eq_example_maximum_likelihood_single_ring}?
To answer this question we can measure the `quality' of the $k_{1}$-image-ring by using the following simple objective-function:
\begin{eqnarray}
\Ccost_{1} & = & \frac{1}{4\fsigma^{2}} \int_{0}^{2\pi}\int_{0}^{2\pi} | \fA^{\signal}(k_{1},\psi) - \fA^{\signal}(k_{1},\psi^{\prime}) |^{2} d\psi d\psi^{\prime} \text{,}
\end{eqnarray}
which, up to a normalization constant, is equivalent to:
\begin{eqnarray}
\Ccost_{1} & = & \frac{1}{4\fsigma^{2}} \int_{0}^{2\pi}\int_{0}^{2\pi} \left\{ \int_{0}^{2\pi} | \rotation_{\gamma}\circ\fA^{\signal}(k_{1},\psi) - \rotation_{\gamma^{\prime}}\circ\fA^{\signal}(k_{1},\psi) |^{2} d\psi \right\} d\gamma d\gamma^{\prime} \text{.}
\end{eqnarray}
The value of $\Ccost_{1}$ is -- up to an additive constant -- equal to the negative-log-probability that a collection of observations of the image-ring $\fA^{\signal}(k_{1},\psi)$ (each sampled at independent and uniformly-distributed values of $\psi$) is equal to another collection of observations of the image-ring $\fA^{\signal}(k_{1},\psi^{\prime})$ (each sampled at independent and uniformly-distributed values of $\psi^{\prime}$).
Equivalently, we can think of $\Ccost_{1}$ as an affine-transformation of the negative-log-probability that two independent randomly-rotated observations of the image-ring $\fA^{\signal}(k_{1},\cdot)$ will equal one another, given random values of the rotation-angles chosen uniformly in $[0,2\pi)$.
If $\Ccost_{1}$ is low, then the $k_{1}$-image-ring contains very little useful information for alignment: a typical randomly-rotated observation of the ring $\fA^{\signal}(k_{1},\cdot)$ can easily be confused for another.
On the other hand, if $\Ccost_{1}$ is high, then the $k_{1}$-image-ring contains useful information: the ring $\fA^{\signal}(k_{1},\cdot)$ will be quite different from most other rotated versions of itself, and the probability of confusing one orientation for another will be quite small.
A similar measure of quality for the $k_{2}$-image-ring is obtained by replacing $k_{1}$ with $k_{2}$:
\begin{eqnarray}
\Ccost_{2} & = & \frac{1}{4\fsigma^{2}} \int_{0}^{2\pi}\int_{0}^{2\pi} \left\{ \int_{0}^{2\pi} | \rotation_{\gamma}\circ\fA^{\signal}(k_{2},\psi) - \rotation_{\gamma^{\prime}}\circ\fA^{\signal}(k_{2},\psi) |^{2} d\psi \right\} d\gamma d\gamma^{\prime} \text{.}
\end{eqnarray}
The ring with the higher quality is the one we should use in \eqref{eq_example_maximum_likelihood_single_ring}.

\subsection{Principal-image-rings}
\label{sec_Principal_image_rings}

Now let's imagine that we are given a choice of not merely measuring one of the image-rings for a specific $k_{r}$, but rather a more general linear-combination of the rescaled image-rings taken across all $k$.
In terms of notation, we'll assume that $\vu=\transpose{[u_{1},u_{2}]}$ is a unit-vector (i.e., $\|\vu\|=1$) and we'll define $[\transpose{\vu}\fA](\psi_{q})$ as:
\begin{eqnarray}
[\transpose{\vu}\fA](\psi_{q}) = \sum_{r=1}^{\rmax} u_{r}\fA(k_{r},\psi_{q})\eta_{r} \text{,}
\end{eqnarray}
where the rescaling-factor $\eta_{r}=\sqrt{w_{r}\dpsi}$ (see \ref{sec_Image_similarity}).
In the two-ring case-study mentioned above, the linear-combination takes the form:
\begin{eqnarray}
[\transpose{\vu}\fA](\psi_{q}) = u_{1}\fA(k_{1},\psi_{q})\eta_{1} + u_{2}\fA(k_{2},\psi_{q})\eta_{2} \text{,}
\end{eqnarray}
with $\eta_{1}=\eta_{2}=\sqrt{\dpsi}$.

In a moment we will choose our $\vu$ to be one of the principal-vectors of an $\rmax\times\rmax$ matrix, and we'll refer to the linear combination $[\transpose{\vu}\fA](\psi_{q})$ as the `principal-image-ring' associated with that principal-vector $\vu$.
To foreshadow this perspective we note that, because the noise is iid across image-rings, the variance of the noise $[\transpose{\vu}\fA^{\noise}](\psi_{q})$ will be equal to $\sum_{r}\vu_{r}^{2}\fsigma^{2}$, which is equal to $\|\vu\|^{2}\fsigma^{2}$.
Because $\vu$ is a unit-vector, this last expression is simply $\fsigma^{2}$, which is the same as the variance of any individual term $\fA^{\noise}(k_{r},\psi_{q})$ for any particular $k_{r}$.
Moreover, if we consider two orthonormal vectors $\vu_{1}$ and $\vu_{2}$, then $[\transpose{\vu_{1}}\fA^{\noise}](\psi_{q})$ and $[\transpose{\vu_{2}}\fA^{\noise}](\psi_{q})$ will be independent random variables, each drawn from $\cN(0,\fsigma^{2})$.

If we were to use the principal-image-rings $[\transpose{\vu}\fA]$ and $[\transpose{\vu}\fB]$ to align the two images, we would consider the approximate maximum-likelihood estimate:
\begin{eqnarray}
\gamma_{\optimal}^{\estimated}(\vu) & = & \argmax_{\gamma} - \sum_{q=0}^{\qmax} | [\transpose{\vu}\fA](\psi_{q}) - [\transpose{\vu}\rotation_{\gamma}\fB](\psi_{q}) |^{2} \text{,}
\end{eqnarray}
where `off-grid' values of $\gamma$ not equal to a multiple of $\dpsi$ are treated via \eqref{eq_fourier_bessel_rotation}.
As the number of angular-quadrature-nodes $\qmax\rightarrow\infty$, this maximum-likelihood-estimate converges to:
\begin{eqnarray}
\gamma_{\optimal}^{\estimated}(\vu) & = & \argmax_{\gamma} - \int_{0}^{2\pi} | \sum_{r} \left[\vu_{r}\fA(k_{r},\psi)\sqrt{w_{r}}\right] - \sum_{r} \left[\vu_{r}\rotation_{\gamma}\fB(k_{r},\psi)\sqrt{w_{r}}\right] |^{2} d\psi \text{.}
\label{eq_example_maximum_likelihood_single_u}
\end{eqnarray}
Note that this is a straightforward generalization of our previous (and more limited) approximation \eqref{eq_example_maximum_likelihood_single_ring}.
Indeed, \eqref{eq_example_maximum_likelihood_single_ring} is obtained by considering \eqref{eq_example_maximum_likelihood_single_u} with $\vu$ chosen to be a column of the identity-matrix.

At this point we remark that the operation-count required to compress any image $\fA$ onto the principal-image-ring $[\transpose{\vu}\fA]$ is $\bigO(\rmax\qmax)$.
This is the same as the operation-count of a single image-image alignment, which is dominated by the operation-count of \eqref{eq_image_innerproduct}-1.
Thus, if the goal is a single image-image alignment, the construction of a principal-image-ring may not seem computationally advantageous.
However, in a typical cryo-EM application, the same image is often aligned many times against many other different images.
In this standard scenario the construction of the principal-image-rings can be performed as a precmputation, with a total operation-count that is negligible when the number of images and targets is large.
See section \ref{sec_Image_alignment_example} for more details.

\subsection{Image objective-function for a single target}
\label{sec_Image_objective_single}

Given the freedom to choose $\vu$, the natural question is: which choice of $\vu$ is the best for alignment?
By generalizing our previous notion of `quality' above, we can try and find the $\vu$ that maximizes:
\begin{eqnarray}
\Ccost(\vu ; \fB) & = & \frac{1}{4\fsigma^{2}} \int_{0}^{2\pi} \int_{0}^{2\pi} \left\{ \sum_{q=0}^{\qmax-1} \left| [\transpose{\vu}\rotation_{\gamma}\fB](\psi_{q}) - [\transpose{\vu}\rotation_{\gamma{\prime}}\fB](\psi_{q}) \right|^{2} \right\} d\gamma d\gamma^{\prime} \text{,}
\label{eq_Image_single_cost}  
\end{eqnarray}
Where we use the term $\fB$ to denote $\fA^{\signal}$.
Similar to before, the value of $\Ccost(\vu ; \fB)$ is an affine-transformation of the negative-log-probability that an observation of the principal-image-ring $[\transpose{\vu}\fB]$ will equal a randomly-rotated observation of the same principal-image-ring.
For the two-ring case-study, choices of $\vu$ for which $\Ccost(\vu ; \fB)$ is low will correspond to linear combinations of the $k_{1}$- and $k_{2}$-image-rings which are not sensitive to the image-orientation, while choices of $\vu$ for which $\Ccost(\vu ; \fB)$ is high will correspond to linear combinations which are highly informative for alignment.
%This observation strongly suggests a choice of $\vu$ which maximizes $\Ccost(\vu ; \fB)$.

While there are many different measures of `quality' on might consider for $\vu$, we find that $\Ccost(\vu ; \fB)$ is particularly convenient because $\Ccost(\vu ; \fB)$ is a quadratic function of $\vu$ that is easy to maximize.
More specifically, we can rewrite $\Ccost(\vu ; \fB)$ as:
\begin{eqnarray}
\Ccost(\vu ; \fB) & = & \sum_{r^{\prime}=1}^{\rmax} \sum_{r=1}^{\rmax} \vu_{r} \cdot \Ckern_{r,r^{\prime}}(\fB) \cdot \vu_{r^{\prime}} \\
              & = & \transpose{\vu}\cdot \Ckern \cdot \vu \text{,}
\end{eqnarray}
where $\Ckern(\fB)$ is a real $\rmax\times\rmax$ symmetric positive-definite matrix.

In this form one can immediately recognize $\Ccost(\vu ; \fB)$ as a rayleigh-quotient of the kernel $\Ckern(\fB)$.
The $\vu$ which maximizes the rayleigh-quotient $\Ccost(\vu ; \fB)$ is equal to the dominant principal-vector (i.e., eigenvector) of the kernel $\Ckern(\fB)$.
Moreover, the dominant principal-value (i.e., eigenvalue) of $\Ckern(\fB)$ will be equal to the quality $\Ccost(\vu ; \fB)$ for that principal-vector.
The next best orthonormal choice of $\vu$ will be the second principal-vector, and so forth, with the orthornormal sequence of principal-vectors corresponding to the principal-vectors (i.e., eigenvectors) of $\Ckern(\fB)$.
%If one were interested, the final principal-vector of \Ckern corresponds to the principal-image-ring with the lowest value of $\Ccost$; this principal-image-ring maximuzes the average log-probability of confusing one orientation for another during alignment.

The matrix entries of $\Ckern(\fB)$ take the form:
\begin{eqnarray}
\Ckern_{r,r^{\prime}}(\fB) & = & \frac{1}{4\fsigma^{2}} \int_{0}^{2\pi} \int_{0}^{2\pi} \sum_{q=0}^{\qmax-1} \ldots \\
 & & \quad \left[ \rotation_{\gamma}\circ\fB(k_{r},\psi_{q})\sqrt{w_{r}} - \rotation_{\gamma^{\prime}}\circ\fB(k_{r},\psi_{q})\sqrt{w_{r}} \right]^{\dagger} \times \ldots \\
 & & \qquad \left[ \rotation_{\gamma}\circ\fB(k_{r^{\prime}},\psi_{q})\sqrt{w_{r^{\prime}}} - \rotation_{\gamma^{\prime}}\circ\fB(k_{r^{\prime}},\psi_{q})\sqrt{w_{r^{\prime}}} \right] \dpsi d\gamma d\gamma^{\prime} \text{.}
\label{eq_Image_single_kern}
\end{eqnarray}
Using Plancherel's theorem in 1-dimension, we note that
\begin{eqnarray}
\int_{0}^{2\pi} \int_{0}^{2\pi} \sum_{q=0}^{\qmax-1} \rotation_{\gamma}\circ\fA(k_{r},\psi_{q}){\dagger} \cdot \rotation_{\gamma}\circ\fB(k_{r^{\prime}},\psi_{q}) \dpsi d\gamma d\gamma^{\prime} = (2\pi)^3 \sum_{q=0}^{\qmax-1} \bA(k_{r},q)^{\dagger} \bB(k_{r^{\prime}},q) \text{,}
\label{eq_image_kernel_ring_same}
\end{eqnarray}
and similarly
\begin{eqnarray}
\int_{0}^{2\pi} \int_{0}^{2\pi} \sum_{q=0}^{\qmax-1} \rotation_{\gamma}\circ\fA(k_{r},\psi_{q}){\dagger} \cdot \rotation_{\gamma^{\prime}}\circ\fB(k_{r^{\prime}},\psi_{q}) \dpsi d\gamma d\gamma^{\prime} = (2\pi)^3  \bA(k_{r},0)^{\dagger} \bB(k_{r^{\prime}},0) \text{.}
\label{eq_image_kernel_ring_diff}
\end{eqnarray}
These two identities allow us to evaluate each matrix entry of $\Ckern_{r,r^{\prime}}(\fB)$ easily:
\begin{eqnarray}
\Ckern_{r,r^{\prime}}(\fB) & = & \frac{(2\pi)^{3}}{4\fsigma^{2}} \sum_{q=1}^{\qmax-1} \bB(k_{r},q)^{\dagger} \bB(k_{r^{\prime}},q) \text{,}
\label{eq_Image_simple_kern}
\end{eqnarray}
where the $q$-sum does not include the $q=0$ term.

Once we have the $\rmax\times\rmax$ matrix $\Ckern(\fB)$, we can take its singular-value-decomposition to define the sequence of principal-vectors $\vu_{1},\ldots,\vu_{\nrank}$ (for some fixed maximum rank $\nrank$).
We can then use these principal-vectors to approximate the inner-product $\cX(\gamma ; \fA,\fB)$.
The steps in this approximation are analogous to \eqref{eq_image_innerproduct}:
\begin{flalign}
\label{eq_image_innerproduct_pm}
& \textbf{[Step 1]} \quad \fourier{\cX}(q) = 2\pi \sum_{r=1}^{\nrank} \left[\transpose{\vu_{r}}\bA\right](q)^{\dagger} \left[\transpose{\vu_{r}}\bB\right](q) \text{,} &&\\\nonumber
& \textbf{[Step 2]} \quad \cX(\gamma_{q^{\prime}}) = \sum_{q=0}^{\qmax-1} \exp\left(-\imunit q \gamma_{q^{\prime}}\right) \fourier{\cX}(q) \text{.} &&
\end{flalign}

One instance of the two-ring case-study described above is illustrated in Fig \ref{pm_fig_image_pair_S2_FIGB}, using two image-rings from the $\fA^{\signal}$ presented in Fig \ref{pm_fig_image_pair_FIGA}.
These two image-rings are chosen so that they are both quite useful for alignment (see dashed- and solid-black lines).
By using the radial-SVD, we can identify a linear-combination $\vu$ of the two rings (shown in red) which is even more useful than either of the image-rings individually.

Another two-ring case-study is shown in Fig \ref{pm_fig_image_pair_S2_FIGC}, using two different image-rings from the same signal.
These two image-rings are chosen so that the first has a very large amplitude (dashed-black), but is not particularly useful for alignment.
The radial-SVD produces a linear-combination $\vu$ which gives the other image-ring (solid-black) a larger weight, even though the amplitude of this second image-ring is smaller than the first.
The resulting linear-combination (red) is the dominant principal-image-ring, which is more useful for alignment than the next principal-image-ring (cyan).

If, instead of limiting ourselves to merely two rings, we consider the entire image shown in Fig \ref{pm_fig_image_pair_FIGA}, we can apply the same technique.
An illustration of the eigenvalues of the kernel $\Ckern(CTF\odot\fS)$ is shown in Fig \ref{pm_fig_image_pair_Sall_FIGD}; one can clearly see that only a handful of principal-vectors are required to capture most of the structure.
By considering only $\nrank$ principal-vectors, we can approximate the inner-product array $\cX(\gamma ; \fA,CTF\odot\fS)$, as described in \eqref{eq_image_innerproduct_pm}.
These approximate inner-product arrays are shown in Fig \ref{pm_fig_image_pair_X_wSM_FIGE}.
Note that the inner-product landscape rapidly converges as $\nrank$ increases.
Note too that, even when $\cX(\gamma ; \fA,CTF\odot\fS)$ is not very accurate (e.g., when $\nrank=2$ or $3$), the optimal alignment-angle (i.e., the $\argmax$ of $\cX$ over $\gamma$) is still quite accurate.
This is because we have designed our objective-function to ignore the magnitude of the different image-rings, instead prioritizing the ability to discriminate between alignment-angles.
Similar results hold if we use the principal-image-rings of $CTF\odot\fS$ to align the images to one another, as shown in Fig \ref{pm_fig_image_pair_X_wSM_FIGF}.

\subsection{Image objective-function for multiple targets}
\label{sec_Image_objective_multiple}

In many cryo-EM applications we are given a large set of $\nimage$ experimental-images.
A common goal is to align these experimental-images to a set of $\ntarget$ target-images.
These targets might be constructed by averaging specially chosen groups of images (to form what are known as `class-averages'), or by projecting a reference-molecule and convolving with the micrograph-specific CTF (to form `templates').
Typically the number of images $\nimage$ is $10^{3}$-$10^{5}$, while the number of targets $\ntarget$ can be $10^{3}$ or more.
In these contexts it it not always necessary to construct principal-image-rings using a separate objective-function for each target.
Instead, one can often make do with a target-averaged objective-function and kernel, such as:

\begin{eqnarray}
\Ccost(\vu) = \frac{1}{\ntarget}\sum_{j=0}^{\ntarget} \Ccost(\vu ; \fB_{j}) \text{ \ \ \ along with \ \ \ } \Ckern & = & \frac{1}{\ntarget}\sum_{j=0}^{\ntarget} \Ckern(\fB_{j}) \text{,}
\label{eq_image_cost_and_kern}
\end{eqnarray}
where each summand involves a single target via \eqref{eq_Image_single_cost} and \eqref{eq_Image_single_kern}, and the index $j$ sums over all the targets $\fB_{j}$.

\subsection{Image alignment example}
\label{sec_Image_alignment_example}

An example of this approach is shown in Fig \ref{pm_fig_M_and_S_FIGG} and Fig \ref{pm_fig_X_2d_cluster_d0_FIGK}.
This example uses the first $\nimage=1024$ images in the {\tt EMPIAR-10005} dataset.
The images span the first $17$ micrographs in the dataset, corresponding to $17$ distinct contrast-transfer-functions (CTFs).
For each of the distinct CTFs, we construct $\ntarget=993$ targets by first projecting the reference-molecule {\tt emd-5778} onto $\ntarget$ distinct viewing-angles, and then correcting each of the resulting templates with the CTF.
For our numerical experiment, we will align the $\nimage$ images to the $\ntarget$ targets.
This numerical experiment simulates the kinds of image-alignment tasks performed in many standard molecular-reconstruction pipelines \cite{EMAN2,Scheres2012b,Kimanius2016,Bell2016}.
A few of the images and targets are shown in Fig \ref{pm_fig_M_and_S_FIGG}.

We first use the strategy of section \ref{sec_Image_innerproducts_in_a_discrete_setting} (taking into account all $\rmax$ radial quadrature-nodes) to calculate the inner-products $\cX(\gamma ; \fA,\fB )$ for each image-target pair $(A,B)$.
We set this full calculation aside as our `ground-truth'.
We then apply our radial-SVD to the same data-set.
Because each CTF corresponds to a distinct set of $\ntarget$ targets, we construct the objective-function and kernel for each CTF separately using \eqref{eq_image_cost_and_kern}.
For each CTF we use the $\nrank$ dominant principal-image-rings to approximate the array of inner-products as described in \eqref{eq_image_innerproduct_pm}.
We'll denote these approximate inner-products via $\cX^{\estimated}(\gamma ; \fA,\fB ; \nrank)$.

When $\nrank=\rmax$, we recover the full calculation (i.e., $\cX^{\estimated}(\gamma ; \fA,\fB ; \rmax)=\cX(\gamma ; \fA,\fB)$ up to machine precision).
When $\nrank<\rmax$ then $\cX^{\estimated}(\gamma ; \fA,\fB ; \rmax)$ is often still quite accurate.
We can compare the values of $\cX^{\estimated}(\gamma ; \fA,\fB ; \rmax)$ to $\cX(\gamma ; \fA,\fB)$ for each $\nrank$ directly, calculating the relative-error in the frobenius-norm over the entire array, shown on the left of Fig \ref{pm_fig_X_2d_cluster_d0_FIGK}.
As one can see, once $\nrank$ is greater than $16$ or so the approximation is quite accurate.

It is important to note that the radial-SVD is often still useful for alignment even when the inner-product landscape is not accurate in terms of the absolute value.
Indeed, we have structured our objective-function in \eqref{eq_Image_single_cost} to ignore the overall magnitude of each image-ring (which only contributes a constant to the inner-product).
Instead, $C$ prioritizes those image-rings which are useful for alignment and discrimination, even if their overall magnitude is not that large.
Thus, the approximate inner-products $\cX^{\estimated}(\gamma ; \fA,\fB ; \nrank)$ can be correlated with the ground-truth, even when $\nrank$ is small.
This correlation (taken across the entire array) is shown in red in the middle of Fig \ref{pm_fig_X_2d_cluster_d0_FIGK}.
Notably, the correlation reaches over $95\%$ for values of $\nrank\sim 6$, even though $\cX^{\estimated}(\gamma ; \fA,\fB ; \nrank)$ is far from $\cX(\gamma ; \fA,\fB)$ in terms of relative-error.

The utility of the radial-SVD is even more apparent when one considers the `backward-error' associated with alignment.
To formalize this notion, let's first fix an image-target pair $(\fA,\fB)$, and then define the `optimal' alignment-angle for that pair to be:
\begin{eqnarray}
\gamma_{\optimal} & = & \argmax_{\gamma} \cX(\gamma) \text{,}
\end{eqnarray}
with the associated optimal inner-product
\begin{eqnarray}
\cX_{\optimal} & = & \cX(\gamma_{\optimal}) \text{,}
\end{eqnarray}
where we have suppressed the arguments $\fA,\fB$ for readability.
Clearly, all the inner-products $\cX(\gamma)$ are less than (or equal to) the optimal $\cX_{\optimal}$; that is, the optimal alignment is at the $100^{\text{th}}$-percentile in the list of true inner-products.

Now let's define the approximate optimal alignment-angle:
\begin{eqnarray}
\gamma_{\optimal}^{\estimated}(\nrank) & = & \argmax_{\gamma} \cX^{\estimated}(\gamma ; \nrank) \text{.}
\end{eqnarray}
To evaluate this approximate optimal alignment-angle in terms of a backwards error, we can first calculate the true inner-product at this approximate alignment-angle:
\begin{eqnarray}
\cX_{\optimal}^{\estimated}(\nrank) & = & \cX(\gamma_{\optimal}^{\estimated}(\nrank)) \text{,}
\end{eqnarray}
and then calculate the fraction $f$ of other true inner-products $\cX(\cdot)$ which are less than (or equal to) $\cX_{\optimal}^{\estimated}$.
In simpler terms, $f$ is the true percentile of the approximated alignment within the list of true inner-products.

Because our radial-SVD emphasizes alignment-accuracy instead of inner-product-magnitude, it is certainly possible for $\cX_{\optimal}^{\estimated}$ to be close to the `top of the list' of true inner-products (i.e., for $f$ to be close to $1$), even if the approximate inner-product $\cX^{\estimated}( \gamma_{\optimal}^{\estimated}(\nrank) ; \nrank)$ is very far from the true inner-product $\cX_{\optimal}^{\estimated}$.
We observed this phenomenon earlier on in Fig \ref{pm_fig_image_pair_X_wSM_FIGE} and Fig \ref{pm_fig_image_pair_X_wSM_FIGF}, and we see it once again in this larger example as well.
The average value of $f$ (taken across image-target pairs) is shown as a function of $\nrank$ in cyan in the middle of Fig \ref{pm_fig_X_2d_cluster_d0_FIGK}.

In terms of operation-count, we expect that it will take $\bigO(\nimage\ntarget\rmax\qmax)$ operations to calculate $\cX(\gamma ; \fA,\fB)$ for all the image-target pairs using the strategy of section \ref{sec_Image_innerproducts_in_a_discrete_setting}.
For our radial-SVD we have the following contributions to the operation-count of $\cX^{\estimated}(\gamma ; \fA,\fB ; \nrank)$:
\begin{enumerate}
\item Calculate $\Ckern$ using \eqref{eq_Image_simple_kern}: $\bigO(\ntarget\qmax\rmax^{2})$
\item Find the SVD of $\Ckern$: $\bigO(\rmax^3)$
\item Calculate the $\left[\transpose{\vu}\fA\right]$ and $\left[\transpose{\vu}\fB\right]$: $\bigO(\nimage\rmax\qmax\nrank + \ntarget\rmax\qmax\nrank)$
\item Calculate $\fourier{\cX}$ using \eqref{eq_image_innerproduct_pm}-1: $\bigO(\nimage\ntarget\nrank\qmax)$
\item Calculate $\cX$ using \eqref{eq_image_innerproduct_pm}-2: $\bigO(\nimage\ntarget\qmax\log(\qmax))$
\end{enumerate}
Note that steps 1-3 are `precomputations' that do not need to be performed for every image-target pair.
Steps 4-5 are `computations' that need to be performed once per image-target pair.
In the limit as $\nimage$ and $\ntarget$ both become very large, we expect the latter two steps to dominate.
In this limit we hope for a speedup of roughly $\rmax/\nrank$, as shown in cyan in the right subplot of Fig \ref{pm_fig_X_2d_cluster_d0_FIGK}.
The actual speedup in runtime (including precomputations), for a straightforward implementation on a dell laptop with an i7 processor, is less dramatic, and is shown in red.
Note that even though steps 1-3 above were repeated for each of the distinct CTFs, the total precomputation-time was still negligible.

\section{Volume alignment}
\label{sec_Volume_alignment}

As mentioned earlier on, we can certainly apply the same techniques to align 3-dimensional volumes rather than 2-dimensional image.
The same general strategy applies, with the only major change being that we can now compress both the radial magnitude $k$ as well as the spherical-harmonic-degree $l$ -- as both are preserved under rotation.
In this section we describe this strategy, highlighting the differences between the 2- and 3-dimensional situations.
For clarity we will typically ignore the constant factors associated with fourier-transforms and integration on the sphere.

\subsection{Volume notation}
\label{sec_Volume_notation}

We use $\vx,\vk \in\Real^{3}$ to represent spatial position and frequency, respectively.
in spherical-coordinates the vector $\vk$ is represented as:
\begin{eqnarray}
\vk &=& k\cdot \hk, \text{\ \ with\ \ } \hk = (\cos \kazimub \sin \kpolara, \sin \kazimub \sin \kpolara , \cos \kpolara) \text{,}
\end{eqnarray}
with polar-angle $\kpolara$ and azimuthal-angle $\kazimub$ representing the unit vector $\hk$ on the surface of the sphere $S^{2}$.

Using a right-handed basis, a rotation about the third axis by angle $\eazimub$ is represented as:
\begin{equation}
\rotation_{\eazimub}^{z} =
\left(
\begin{array}{ccc}
+\cos\eazimub & -\sin\eazimub & 0 \\
+\sin\eazimub & +\cos\eazimub & 0 \\
0            & 0            & 1
\end{array}
\right) \text{,}
\end{equation}
and a rotation about the second axis by angle $\epolara$ is represented as:
\begin{equation}
\rotation_{\epolara}^{y} =
\left(
\begin{array}{ccc}
+\cos\epolara & 0 & +\sin\epolara \\
0            & 1 & 0            \\
-\sin\epolara & 0 & +\cos\epolara 
\end{array}
\right) \text{.}
\end{equation}

A rotation $\rotation_{\tau}$ of a vector $\vk\in\Real^{3}$ can be represented by the vector of euler-angles $\tau = (\egammaz,\epolara,\eazimub)$:
\begin{equation}
\rotation_{\tau}\cdot\vk = \rotation_{\eazimub}^{z}\circ \rotation_{\epolara}^{y}\circ \rotation_{\egammaz}^{z}\cdot \vk \text{.}
\end{equation}

We represent any given volume as a function $A\in\Leb^{2}(\Real^{3})$, with values corresponding to the volume intensity at each location $\vx\in\Omega_{1}$.
We'll refer to $\fA(\vk)$ in spherical-coordinates as $\fA(k, \hk)$.
With this notation, each $\fA(k,\cdot)$ corresponds to a `shell' in frequency-space with radius $k$.
The rotation of any volume $\rotation_{\tau}\fA(k,\hk)$ corresponds to the function $\fA(k,\rotation_{\tau}^{-1}\hk)$.

\subsection{Volume inner-products in a continuous setting: the spherical-harmonic basis}
\label{sec_Volume_innerproducts_in_a_continuous_setting}

Using the notation above, we can represent a volume $\fA(k,\hk)$ as:
\begin{equation}
\fA(k,\hk) = \sum_{l=0}^{+\infty}\sum_{m=-l}^{m=+l} \fA_{l}^{m}(k) Y_{l}^{m}(\hk) \text{,}
\end{equation}
where $Y_{l}^{m}(\hk)$ represents the spherical-harmonic of degree-$l$ and degree-$m$:
\begin{equation}
Y_{l}^{m}(\kpolara,\kazimub) =  Z_{l}^{m} \euler^{+\imunit\kazimub}P_{l}^{m}(\cos\kpolara) \text{,}
\end{equation}
with $P_{l}^{m}$ representing the (unnormalized) associated Legendre polynomial, and $Z_{l}^{m}$ the normalization-constant:
\begin{equation}
Z_{l}^{m} = \sqrt{\frac{2l+1}{4\pi}\times\frac{\left(l-|m|\right)!}{\left(l+|m|\right)!}} \text{.}
\end{equation}
The coefficients $\fA_{l}^{m}(k)$ define the spherical-harmonic expansion of the $k$-shell $\fA(k,\cdot)$.

Using this spherical-harmonic basis allows us to efficiently apply rotations.
For example, given a rotation $\tau=(\egammaz,\epolara,\eazimub)$, we can represent the rotated volume $\fB:=\rotation_{\tau}\fA$ as:
\begin{equation}
\fB_{l}^{m_{1}} = \sum_{m_{2}=-l}^{m_{2}=+l} \euler^{-\imunit m_{1} \eazimub}d_{m_{1},m_{2}}^{l}(\epolara) \euler^{-\imunit m_{2} \egammaz} \fA_{l}^{m_{2}} \text{,}
\end{equation}
where $d_{m_{1},m_{2}}^{l}(\epolara)$ represents the degree-$l$ wigner-d matrix associated with the interior euler-angle $\epolara$.

Given any two volumes $\fA$ and $\fB$, the inner-product
\begin{eqnarray}
\cX(\tau ; \fA,\fB) & := & \langle \fA , \rotation_{\tau}\fB \rangle \\
                    & =  & \iiint \fA(k,\hk)^{\dagger} \fB(k,\rotation_{\tau}^{-1} \hk) k^{2}dkd\hk
\end{eqnarray}
can be rewritten (up to a constant factor) as:
\begin{eqnarray}
\cX(\tau ; \fA,\fB) & = & \int_{k=0}^{+\infty} \sum_{l=0}^{+\infty}\sum_{m_{1}=-l}^{m_{1}=+l}\sum_{m_{2}=-l}^{m_{2}=+l} \euler^{-\imunit m_{1} \eazimub}d_{m_{1},m_{2}}^{l}(\epolara)\euler^{-\imunit m_{2} \egammaz} \fA_{l}^{m_{1}}(k)^{\dagger} \fB_{l}^{m_{2}}(k) k^{2}dk \\
                    & = & \sum_{m_{1}=-\infty}^{m_{1}=+\infty}\sum_{m_{2}=-\infty}^{m_{2}=+\infty} \euler^{-\imunit m_{1} \eazimub} \euler^{-\imunit m_{2} \egammaz} \left[ \sum_{l=0}^{+\infty} d_{m_{1},m_{2}}^{l}(\epolara) \int_{k=0}^{+\infty} \fA_{l}^{m_{1}}(k)^{\dagger} \fB_{l}^{m_{2}}(k) k^{2}dk \right] \text{.}
\end{eqnarray}
This last expression can be interpreted as a relationship between the desired inner-products $\cX(\tau)$ and the 2-dimensional fourier-transform of the term in brackets on the right-hand-side:
\begin{eqnarray}
\fourier{\cX}(m_{1},m_{2} ; \epolara ; \fA,\fB) = \sum_{l=0}^{+\infty} d_{m_{1},m_{2}}^{l}(\epolara) \left[ \int_{k=0}^{+\infty}\fA_{l}^{m_{1}}(k)^{\dagger} \fB_{l}^{m_{2}}(k) k^{2}dk \right] \text{.}
\end{eqnarray}
Note that the right-hand-side involves accumulating information over both the radius $k$ as well as the spherical-harmonic-degree $l$.
We will aim to compress both of these later on.

\subsection{Volume discretization}
\label{sec_Volume_discretization}

Similar to our discretization of images, we assume that the volume $A(\vx)$ is supported in $\Omega_{1}$, and that most of the relevant frequency-content is contained in $\Omega_{\kmax}$.

With these assumptions, we can disretize the radial component of $\Omega_{\kmax}\in\Real^{3}$ using a Gauss-Jacobi quadrature for $k$ built with a weight-function corresponding to a radial-weighting of $k^{2}dk$; once again the number of radial quadrature-nodes $\rmax$ will be $\bigO(\kmax)$.
Each of the shells $\fA(k_{r},\cdot)$ can be accurately described using spherical-harmonics with $l\leq\bigO(k_{r})$.
Thus, the number of spherical-harmonic coefficients required for each shell $\fA(k_{r},\cdot)$ is $\bigO(k_{r}^{2})$.
The total number of spherical-harmonic coefficients required to approximate $\fA(k,\hk)$ over $\Omega_{\kmax}$ is $\bigO(\kmax^{3})$, with a maximum degree of $\lmax=\bigO(\kmax)$.
The associated maximum order will be $\mmax=1+2\lmax$, which is also $\bigO(\kmax)$.
For brevity we will treat the order indices $m_{1}$ and $m_{2}$ periodically in the interval $[-\lmax,\ldots,+\lmax]$ (so, for example, the $m$-value of $\lmax+1$ corresponds to the $m$-value of $-\lmax$); coefficients $\fA_{l}^{m}$ with $l<|m|$ will be identically $0$.
%If necessary, we typically discretize $\hk$ for each shell using legendre-nodes for $\cos\kpolara$ to define latitudinal rings (and associated latitudinal weights), along with equispaced nodes for $\kazimub$ on each latitudinal ring.

\subsection{Volume inner-products in a discrete setting}
\label{sec_Volume_innerproducts_in_a_discrete_setting}

Using the formulae of section \ref{sec_Volume_innerproducts_in_a_continuous_setting}, we can calculate the inner-products $\cX$ for any interior-angle $\epolara$ across a range of azimuthal-angles $\eazimub_{m_{1}'}=2\pi m_{1}'/\mmax$ and in-plane-angles $\egammaz_{m_{2}'}=2\pi m_{2}'/\mmax$. The calculation can be summarized as:

\begin{equation}
\cX(\tau_{m_{1}',m_{2}'} ; \fA,\fB) = \sum_{m_{1}=0}^{\mmax-1} \sum_{m_{2}=0}^{\mmax-1} \exp(-\imunit m_{1} \eazimub_{m_{1}'}) \exp(-\imunit m_{2} \egammaz_{m_{2}'}) \sum_{l=0}^{\lmax} d_{m_{1},m_{2}}^{l}(\epolara) \sum_{r=1}^{\rmax} w_{r} \fA_{l}^{m_{1}}(k_{r})^{\dagger}\fB_{l}^{m_{2}}(k_{r}) \text{,}
\end{equation}
where the euler-angle $\tau_{m_{1}',m_{2}'} = (\egammaz_{m_{2}'},\epolara,\eazimub_{m_{1}'})$.

This calculation can be broken into the following three steps:
\begin{flalign}
\label{eq_volume_innerproduct}
& \textbf{[Step 1]} \quad \tilde{\cX}(m_{1},m_{2} ; l ; \fA,\fB) = \sum_{r=1}^{\rmax} w_{r} \fA_{l}^{m_{1}}(k_{r})^{\dagger}\fB_{l}^{m_{2}}(k_{r}) \text{,} && \\\nonumber
& \textbf{[Step 2]} \quad \fourier{\cX}(m_{1},m_{2} ; \epolara ; \fA,\fB) = \sum_{l=0}^{\lmax} d_{m_{1},m_{2}}^{l}(\epolara) \tilde{\cX}(m_{1},m_{2} ; l ; \fA,\fB) \text{,} && \\\nonumber
& \textbf{[Step 3]} \quad \cX(\tau_{m_{1}',m_{2}'} ; \fA,\fB) = \sum_{m_{1}=0}^{\mmax-1} \sum_{m_{2}=0}^{\mmax-1} \exp(-\imunit m_{1} \eazimub_{m_{1}'}) \exp(-\imunit m_{2} \egammaz_{m_{2}'}) \fourier{\cX}(m_{1},m_{2} ; \epolara ; \fA,\fB) \text{.} &&
\end{flalign}

The first step combines information from different $k$-shells, requiring $\bigO(\rmax\lmax\mmax^{2})$ operations, but only needs to be performed once per volume-target pair; the results can be reused for different values of $\epolara$.
The second step combines information across different degrees, requiring $\bigO(\lmax\mmax^{2})$ operations for each value of $\epolara$.
The third step can be evaluated using a 2-dimensional fast-fourier-transform of size $\mmax$, requiring $\bigO(\mmax^{2}\log(\mmax))$ operations for each $\epolara$.
Note that the first and second steps dominate the operation-count, as they are more expensive than the third.
In order to resolve the landscape of inner-products over all $\tau$, it is typically necessary to repeat the second and third steps for a grid of $\mmax$ different $\epolara$-values in the interval $[-\pi,+\pi]$.
Consequently, the total operation-count (using all $\rmax$ radial quadrature-nodes and all $\lmax$ degrees) is $\bigO(\kmax^{4})$.

\subsection{Principal-volume-shells}
\label{sec_Principal_volume_shells}

The calculation described in \eqref{eq_volume_innerproduct} involves a $\epolara$-independent sum over the radial quadrature-nodes $k_{r}$, followed by a $\epolara$-dependent sum over the spherical-harmonic-degree $l$.
In order to accelerate this computation, we will need to compress both sums.

We'll start with the radial-compression.
In much the same manner as section \ref{sec_Principal_image_rings}, we can define the principal-volume-shells $[\transpose{\vu}\fA](\hk)$ as:
\begin{eqnarray}
[\transpose{\vu}\fA](\hk) = \sum_{r=1}^{\rmax} u_{r}\fA(k_{r},\hk)\eta_{r} \text{.}
\end{eqnarray}
The $[\transpose{\vu}\fA]$ can be represented using spherical-harmonic coefficients as well, with:
\begin{eqnarray}
[\transpose{\vu}\fA]_{l}^{m} = \sum_{r=1}^{\rmax} u_{r}\fA_{l}^{m}(k_{r})\eta_{r} \text{.}
\end{eqnarray}
In these expressions the rescaling-factor $\eta_{r}$ should be proportional to the square-root of the radial quadrature-weights $w_{r}$.
In the 2-dimensional case of images these weights $w_{r}$ accounted for the radial-weighting of $kdk$.
This time, in the 3-dimensional case of volumes, these weights $w_{r}$ will account for the radial-weighting of $k^{2}dk$ associated with volumetric noise in $\Real^{3}$.

\subsection{Volume objective-function for a single target}
\label{sec_Volume_objective_single}

Given a target-volume of $\fB:=\fA^{\signal}$, we can choose the radial principal-vector $\vu$ to maximize the simple objective-function:
\begin{equation}
\Ccost(\vu ; \fB) = \int_{SO3} \int_{SO3} \sum_{l=0}^{\lmax}\sum_{m=0}^{\mmax}\left| \left[\transpose{\vu}\rotation_{\tau}\fB\right]_{l}^{m} - \left[\transpose{\vu}\rotation_{\tau'}\fB\right]_{l}^{m} \right|^{2} d\tau d\tau' \text{,}
\label{eq_Volume_single_radial_cost}
\end{equation}
where we integrate over the rotations $\tau$ and $\tau'$ using the uniform measure on the group of rotations $SO3$.
Just as before, this objective-function is an affine-transformation of the negative-log-probability that an observation of the principal-volume-shell $[\transpose{\vu}\fB]$ will equal a randomly rotated version of that same shell.

Dropping constant factors, the associated radial kernel $\Ckern(\fB)$ takes the form:
\begin{eqnarray}
\Ckern_{r,r'}(\fB) & = & \int_{SO3} \int_{SO3} \sum_{l=0}^{\lmax}\sum_{m=0}^{\mmax-1} \ldots \\
 & & \quad \left[ \left[\rotation_{\tau}\fB\right]_{l}^{m}(k_{r})\sqrt{w_{r}} - \left[\rotation_{\tau'}\fB\right]_{l}^{m}(k_{r})\sqrt{w_{r}} \right]^{\dagger} \times \ldots \\
 & & \qquad \left[ \left[\rotation_{\tau}\fB\right]_{l}^{m}(k_{r'})\sqrt{w_{r'}} - \left[\rotation_{\tau'}\fB\right]_{l}^{m}(k_{r'})\sqrt{w_{r'}} \right] d\tau d\tau' \text{.}
\label{eq_Volume_single_radial_kern}
\end{eqnarray}
Using the properties of spherical harmonics, as well as the orthogonality relations of the wigner-d matrix \cite{BL81}, we can simplify the radial kernel to:
\begin{eqnarray}
\Ckern_{r,r'}(\fB) = \sum_{l=1}^{\lmax}\sum_{m=0}^{\mmax-1} \fB_{l}^{m}(k_{r})^{\dagger}\fB_{l}^{m}(k_{r'}) \text{,}
\label{eq_Volume_simple_radial_kern}
\end{eqnarray}
where the $l$-sum does not include the $l=0$ term.
As in the 2-dimensional case, we'll use the singular-value-decomposition of the radial kernel $\Ckern$ to define the radial principal-vectors $\vu_{1},\ldots,\vu_{\rmax}$, which will be used below to construct principal-volume-shells.

An important distinction between the 2- and 3-dimensional cases is the overall operation-count of compressing the objects.
In 2-dimensions, the formation of a single principal-image-ring required $\bigO(\rmax\qmax)$ operations, which was equivalent to the operation-count of calculating the inner-products $\cX(\gamma ; \fA,\fB)$ for a single image-target pair using all the $\rmax$ radial quadrature-nodes.
Forming all the principal-image-rings for a single image required $\bigO(\kmax^{3})$ operations, as did finding all the entries of $\Ckern$; each of these tasks required an order-of-magnitude more operations than a single image-target alignment.
Consequently, the radial-SVD was only computationally advantageous in 2-dimensions when aligning multiple images to {\em multiple} targets.

In 3-dimensions, by contrast, the formation of a single principal-volume-shell requires $\bigO(\rmax\lmax\mmax)$ operations, and all the principal-volume-shells for a given volume can be formed in $\bigO(\kmax^{4})$.
This is comparable to the $\bigO(\kmax^{4})$ required to calculate the inner-products $\cX(\tau ; \fA,\fB)$ using all $\rmax$ radial quadrature-nodes and all $\lmax$ degrees for a single volume-target pair.
Moreover, the entries of $\Ckern$ can also be calculated in $\bigO(\rmax^{2}\lmax\mmax)$ operations, which is again comparable to the $\bigO(\kmax^{4})$ operation-count required for a single volume-target alignment.
Thus, as we'll see below, the radial-SVD will be computationally advantageous even when aligning multiple volumes to a {\em single} target; we don't necessarily need multiple volumes {\em and} multiple targets.

\subsection{Principal-volume-degrees}
\label{sec_Principal_volume_degrees}

Now we move on to compressing the spherical-harmonic-degree $l$.
The basic idea is straightforward: the various degrees are each preserved under rotation, yet not all of them are equally important.
Unsurprisingly, we search for linear-combinations of degrees which are useful for alignment.

In terms of notation, we can use any orthonormal vector $\vv\in\Real^{1+\lmax}$ as a weight to accumulate the entries of any array $\vg$ indexed by $l$:
\begin{eqnarray}
[\transpose{\vv}\vg] = \sum_{l=0}^{\lmax} v_{l} g_{l} \text{.}
\label{eq_principal_degree_array}
\end{eqnarray}
When applied to the wigner-d-matrix we get:
\begin{eqnarray}
[\transpose{\vv}d]_{m_{1},m_{2}}(\epolara) = \sum_{l=0}^{\lmax} v_{l} d_{m_{1},m_{2}}^{l}(\epolara) \text{.}
\label{eq_principal_degree_wigner}
\end{eqnarray}
When applied to the spherical-harmonic coefficients of a volume we get:
\begin{eqnarray}
[\transpose{\vv}\fB]^{m}(k_{r}) = \sum_{l=0}^{\lmax} v_{l} \fB_{l}^{m}(k_{r}) \text{.}
\label{eq_principal_degree_spharm}
\end{eqnarray}
Because the different degrees $l$ of $\fB_{l}^{m}(k_{r})$ will typically be associated with the same variance (for any specific value of $|m|\leq l$), there is no need for an $l$-dependent rescaling-factor.

\subsection{Volume objective-function for degrees}
\label{sec_Volume_objective_degree}

To choose the $\vv$ for our principal-volume-degrees, we'll construct an objective-function very similar to the $\Ccost(\vu)$ used in \eqref{eq_Volume_single_radial_cost}.
In this case we choose $\vv\in\Real^{1+\lmax}$ to maximize:
\begin{equation}
\Dcost(\vv ; \fB) = \int_{SO3} \int_{SO3} \sum_{r=1}^{\rmax} \sum_{m=0}^{\mmax}\left| \left[\transpose{\vv}\rotation_{\tau}\fB\right]^{m}(k_{r}) - \left[\transpose{\vv}\rotation_{\tau'}\fB\right]^{m}(k_{r}) \right|^{2} w_{r} d\tau d\tau' \text{.}
\label{eq_Volume_single_degree_cost}
\end{equation}
Just as before, this is an affine-transformation of the negative-log-probability that an observation of the principal-volume-degree $[\transpose{\vv}\fB]$ will equal a randomly rotated version of that same degree.

The associated degree-wise kernel $\Dkern(\fB)$ is (up to constant factors):
\begin{eqnarray}
\Dkern_{l,l'}(\fB) = \delta_{0l}\delta_{0l'}\sum_{r=1}^{\rmax}\sum_{m=0}^{\mmax-1} \fB_{l}^{m}(k_{r})^{\dagger}\fB_{l'}^{m}(k_{r}) w_{r}\text{,}
\label{eq_Volume_simple_degree_kern}
\end{eqnarray}
where the kronecker-$\delta$ prefactors ensure that $\Dkern_{l,l'}$ is $0$ whenever either $l$ or $l'$ is $0$.
We'll choose our $\vv_{1},\ldots,\vv_{\lmax}$ to be the principal-vectors of the degree-wise kernel $\Dkern$.
Below we'll use these to compress various 3-dimensional arrays onto their principal-volume-degrees.

\subsection{Volume alignment using principal-volumes}
\label{sec_Volume_alignment_using_principal_volumes}

Putting together the components above, we can easily approximate the inner-products $\cX(\tau ; \fA,\fB)$.
We start by fixing $\nrank_{\Ccost}$ and $\nrank_{\Dcost}$.
These will be the number of radial- and degree-wise principal-vectors used in the approximation.
The steps in the approximation are largely analogous to those of \eqref{eq_volume_innerproduct}:
\begin{flalign}
\label{eq_volume_innerproduct_pm}
& \textbf{[Step 1]} \quad \tilde{\cX}(m_{1},m_{2} ; l ; \fA,\fB) = \sum_{r=1}^{\nrank_{\Ccost}} \left\{\left[\transpose{\vu_{r}}\fA\right]_{l}^{m_{1}}\right\}^{\dagger}\left[\transpose{\vu_{r}}\fB\right]_{l}^{m_{2}} \text{,} && \\\nonumber
& \textbf{[Step 1b]} \quad \text{form each\ } \left[\transpose{\vv_{l}} \tilde{\cX}\right](m_{1},m_{2} ; \fA,\fB) \quad \forall\ \ l\in\{1,\ldots,\nrank_{\Dcost}\} \text{,} && \\\nonumber
& \textbf{[Step 2]} \quad \fourier{\cX}(m_{1},m_{2} ; \epolara ; \fA,\fB) = \sum_{l=1}^{\nrank_{\Dcost}} \left[\transpose{\vv_{l}}d\right]_{m_{1},m_{2}}(\epolara) \left[\transpose{\vv}_{l}\tilde{\cX}\right](m_{1},m_{2} ; \fA,\fB) \text{,} && \\\nonumber
& \textbf{[Step 3]} \quad \cX(\tau_{m_{1}',m_{2}'} ; \fA,\fB) = \sum_{m_{1}=0}^{\mmax-1} \sum_{m_{2}=0}^{\mmax-1} \exp(-\imunit m_{1} \eazimub_{m_{1}'}) \exp(-\imunit m_{2} \egammaz_{m_{2}'}) \fourier{\cX}(m_{1},m_{2} ; \epolara ; \fA,\fB) \text{.} &&
\end{flalign}

An additional step has been added in betwen steps one and two.
This additional step-1b uses the $\vv_{l}$ to form the principal-degrees of the array $\tilde{\cX}$.
Step-1b needs to be repeated for each volume-target pair, but the results can be used for all values of $\epolara$.

When aligning multiple volumes $\fA$ to a single target $\fB$ (using an array of $\mmax$ $\epolara$-values for each volume-target pair), we have the following contributions to the operation-count:
\begin{enumerate}
\item Calculate $\Ckern(\fB)$ using \eqref{eq_Volume_single_radial_cost}: $\bigO(\lmax\mmax\rmax^{2})$
\item Find the SVD of $\Ckern(\fB)$: $\bigO(\rmax^{3})$
\item Calculate $\Dkern(\fB)$ using \eqref{eq_Volume_single_degree_cost}: $\bigO(\rmax\mmax\lmax^{2})$
\item Find the SVD of $\Dkern(\fB)$: $\bigO(\lmax^{3})$
\item Calculate the $\left[\transpose{\vv}d\right]$: $\bigO(\mmax^{3}\nrank_{\Dcost})$
\item Calculate the $\left[\transpose{\vu}\fB\right]$: $\bigO(\rmax\mmax\lmax\nrank_{\Ccost})$
\item Calculate the $\left[\transpose{\vu}\fA\right]$: $\bigO(\nimage\rmax\mmax\lmax\nrank_{\Ccost})$
\item Calculate $\tilde{\cX}$ using \eqref{eq_volume_innerproduct_pm}-1: $\bigO(\nimage\nrank_{\Ccost}\lmax\mmax^{2})$
\item Calculate the $\left[\transpose{\vv}\tilde{\cX}\right]$ using \eqref{eq_volume_innerproduct_pm}-1b: $\bigO(\nimage\nrank_{\Dcost}\lmax\mmax^{2})$
\item Calculate $\fourier{\cX}$ using \eqref{eq_volume_innerproduct_pm}-2: $\bigO(\nimage\nrank_{\Dcost}\mmax^{3})$
\item Calculate $\cX$ using \eqref{eq_volume_innerproduct_pm}-3: $\bigO(\nimage\mmax^{3}\log(\mmax))$
\end{enumerate}
Note that steps 1-6 are `precomputations' that do not need to be performed for every volume-target pair.
Steps 7-9 are computations that need to be performed once per volume-target pair, but are independent of $\epolara$.
Steps 10-11 are computations that involve the different values of $\epolara$.

When the number of volumes $\nimage$ is large, we expect the total operation-count to scale as $\bigO(\nimage\kmax^{3}\nrank)$, where $\nrank=\max(\nrank_{\Ccost},\nrank_{\Dcost})$, as opposed to the $\bigO(\nimage\kmax^{4})$ operation-count required by \eqref{eq_volume_innerproduct}.
Consequently, as $\nimage$ becomes large we expect a speedup of roughly $\kmax/\nrank$.

\subsection{Volume alignment example}
\label{sec_Volume_alignment_example}

\begin{figure}
\centering
\includegraphics[width=7in]{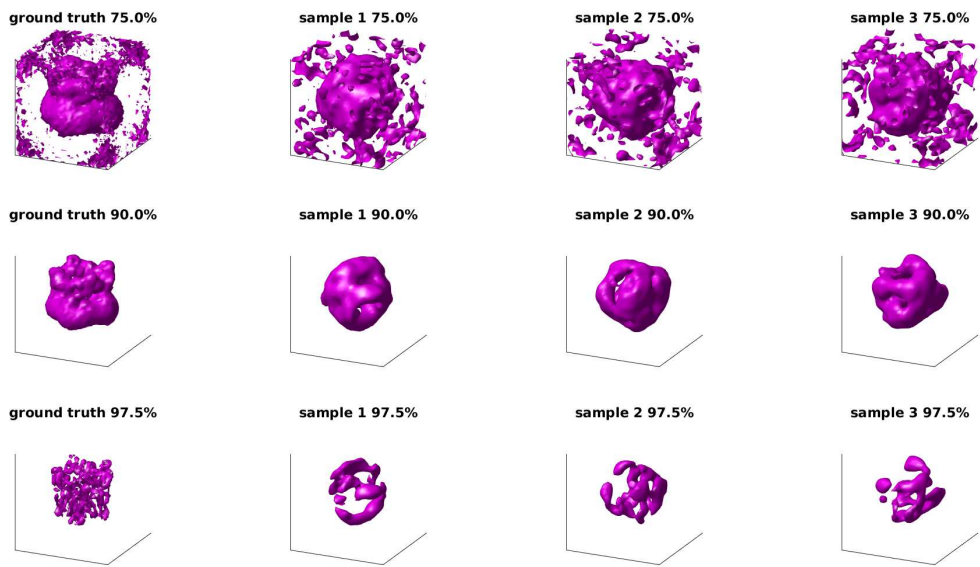}
\caption{\label{pm_vfig_volume_FIGA}
Here we show an array of subplots illustrating different reconstructions of the TRPV1-molecule.
In the leftmost column we show the reference-molecule {\tt emd-5778}.
In the other three columns we show three other volumes, each created using a de-novo reconstruction applied to a small batch of images taken from {\tt EMPIAR-10005}.
Each row corresponds to a different percentile used for the isosurface (i.e., a different level-set of the volume).
}
\end{figure}

Here we present an example illustrating our approach within the context of cryo-EM molecular reconstruction.
For this example we align $\nimage=96$ different volumes to a single target.
Each of the volumes will correspond to an approximation of the TRPV1-molecule produced using a `de-novo' reconstruction-pipeline applied to a subset of $1024$ images taken from the {\tt EMPIAR-10005} dataset.
The target will correspond to the reference-molecule from {\tt emd-5778}.
This numerical experiment simulates the kinds of crossvalidation and/or bootstrapping that are often applied to determine the quality of a dataset and/or a reference-molecule \cite{Henderson2012,Penczek2014,Heymann2015,Rosenthal2015}.
The target is shown in Fig \ref{pm_vfig_volume_FIGA}, alongside a few of the de-novo reconstructions.

\begin{figure}
\centering
\includegraphics[width=6in]{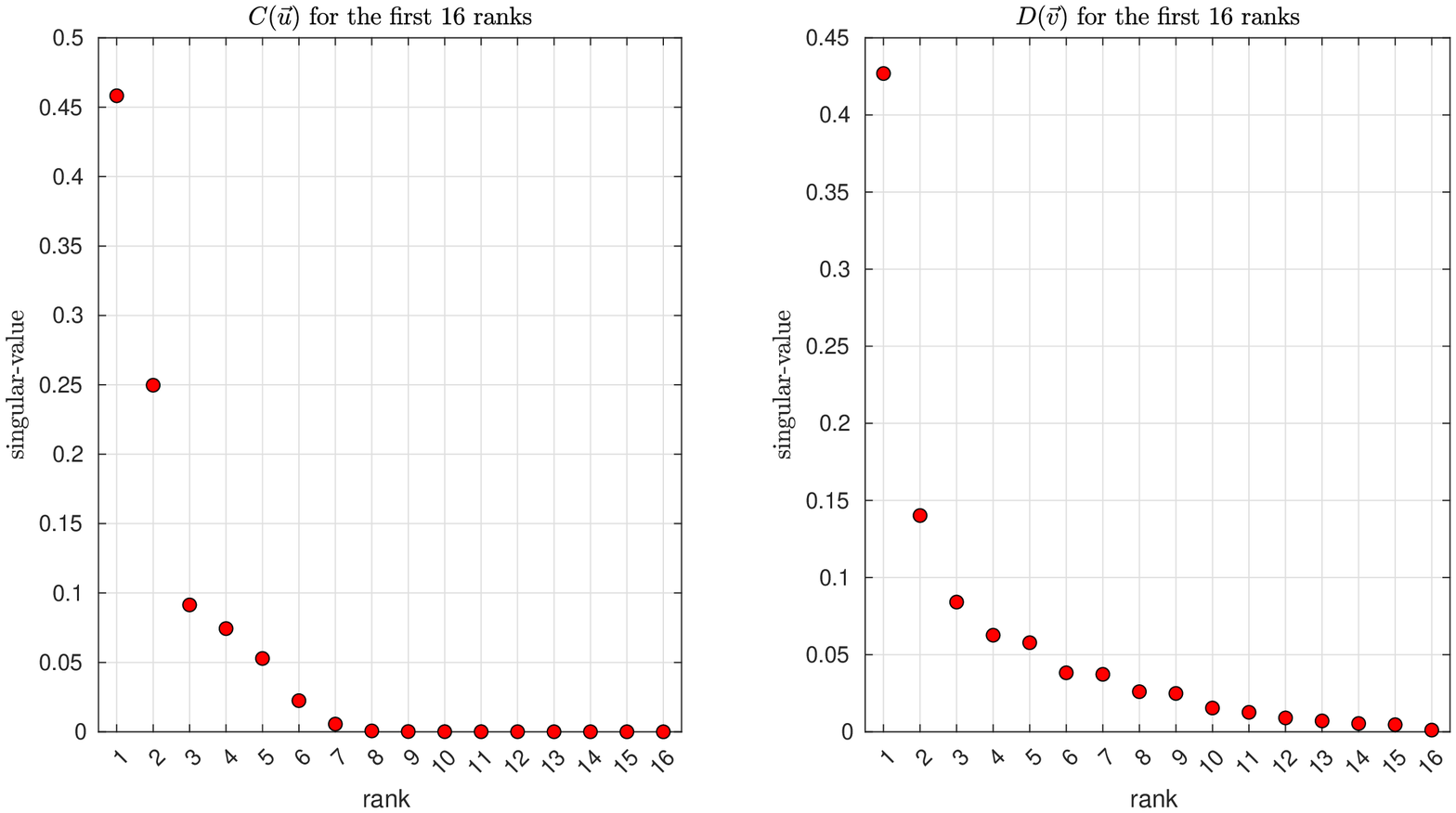}
\caption{\label{pm_vfig_spectrum_FIGD}
Here we show the eigenvalues of $\Ckern$ and $\Dkern$ for the objective-functions $\Ccost(\vu)$ and $\Dcost(\vv)$, respectively.
}
\end{figure}

\begin{figure}
\centering
\includegraphics[width=7in]{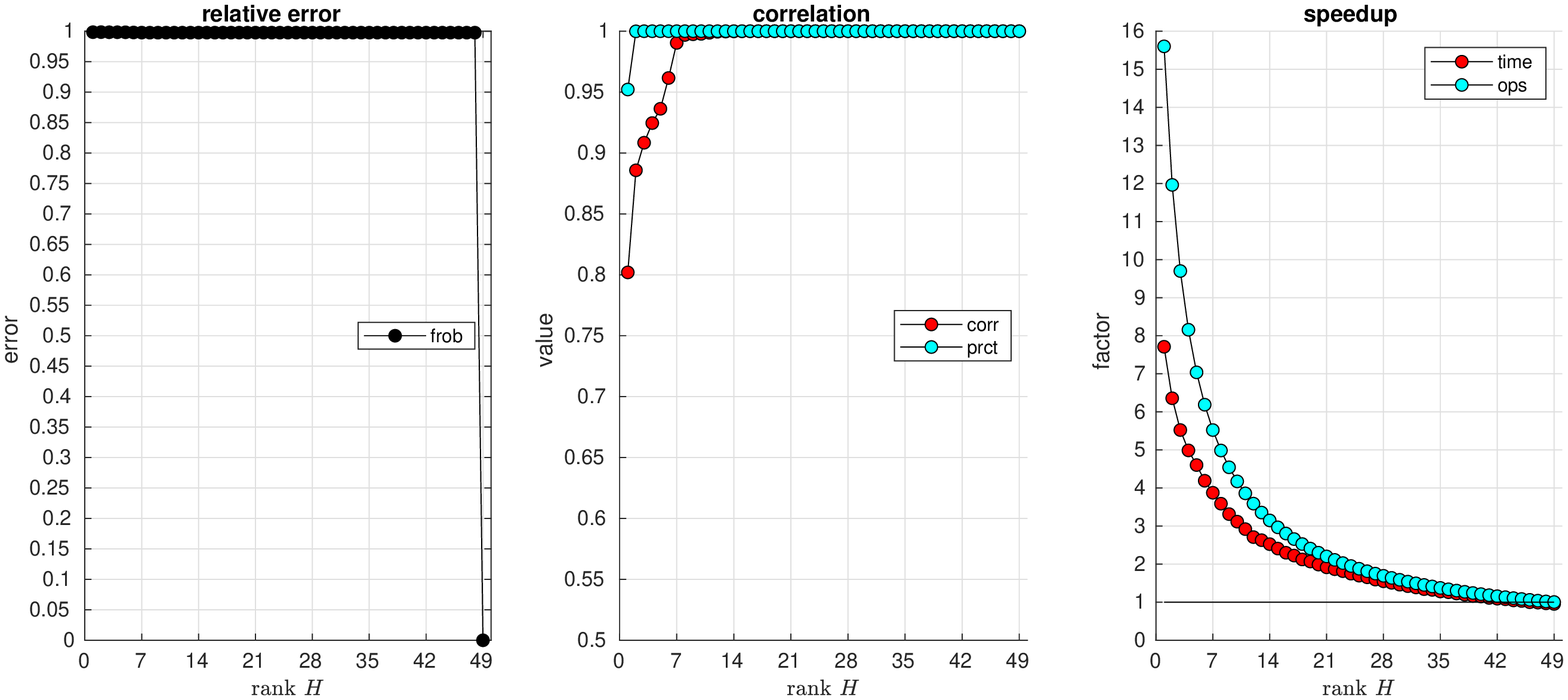}
\caption{\label{pm_vfig_timing_FIGB}
Here we show the results after applying our radial- and degree-SVD to align $\nimage=96$ different volumes to a single target.
Each of the volumes corresponds to a de-novo molecular reconstruction of the TRPV1-molecule using a subset of the {\tt EMPIAR-10005} dataset; three of these are shown on the right of Fig \ref{pm_vfig_volume_FIGA}.
The target is the single reference-molecule shown on the left of Fig \ref{pm_vfig_volume_FIGA}.
We use the strategy of section \ref{sec_Volume_innerproducts_in_a_discrete_setting} to calculate $\cX(\tau ; \fA,\fB )$ for each image-target pair $(A,B)$.
We approximate this inner-product array by using the radial- and degree-SVD, retaining $\nrank$ principal-vectors of $\Ckern$ and $\nrank$ principal-vectors of $\Dkern$, as described in the main text.
For each $\nrank$, we calculate $\cX^{\estimated}(\tau ; \fA,\fB ; \nrank)$, and compare the results to the full calculation (involving all $\rmax$ radial quadrature-nodes and all $\lmax$ degrees).
The relative error (using the frobenius-norm) between these two arrays is shown in black on the left.
Because our principal-volumes are chosen to ignore the rotationally-invariant contributions of the $l=0$ orders, the approximate inner-products $\cX^{\estimated}$ are very far from the true inner-products $\cX$ in terms of absolute-value.
Nevertheless, the $\cX^{\estimated}$ still capture the features of the inner-product-landscape, as reflected in the middle subplot (see main text).
In the limit as the number of image-target pairs becomes very large, we expect a speedup of roughly $\kmax/\nrank$, as shown in cyan in the right subplot.
The actual speedup in total runtime (which includes the necessary precomputations), for an implementation on a dell laptop with an i7 processor, is shown in red.
}
\end{figure}

We apply the strategy of section \ref{sec_Volume_alignment_using_principal_volumes}, using a straightforward implementation with the same choice of $\nrank$ for both $\nrank_{\Ccost}$ and $\nrank_{\Dcost}$.
The results are shown in Fig \ref{pm_vfig_spectrum_FIGD} and Fig \ref{pm_vfig_timing_FIGB}.
Note that the frobenius-norm error associated with our approximation is enormous.
This is because our radial- and degree-wise principal-volumes are chosen to focus on alignment, and discard the rotationally-invariant (but large in magnitude) terms associated with the $k=0$ shells and the $l=0$ degrees.
Nevertheless, the correlation between our approximation and the full calculation is still quite high, shown in red in the middle of Fig \ref{pm_vfig_timing_FIGB}.
Moreover, the `backwards-error' is quite low, as indicated by the high fraction $f$ shown in cyan in the middle of Fig \ref{pm_vfig_timing_FIGB}.
This $f$ is calculated in the same manner as described in section \ref{sec_Image_alignment_example}.
Just as in the case with images, the speedup for our implementation (red) does not quite reach the ideal speedup associated with the operation-count (cyan), shown on the right of Fig \ref{pm_vfig_timing_FIGB}.

\section{Discussion}
\label{sec_Discussion}

For image-alignment, the success of the radial-SVD depends critically on the decay of the spectrum of $\Ckern$.
If the spectrum of $\Ckern(\fB)$ decays slowly -- which will be the case if $\Ckern(\fB)$ is built using a very noisy $\fB$ -- then $\nrank$ will need to be close to $\rmax$ in order to maintain accuracy, and the radial-SVD will not be as computationally advantageous.
On the other hand, if the $\fB$ used to build $\Ckern$ is smooth -- as is often the case in cryo-EM -- then the spectrum of $\Ckern$ will decay relatively quickly, and $\nrank$ can be significantly smaller than $\rmax$ while still maintaining accuracy.
As we demonstrate in section \ref{sec_Image_alignment_example}, the radial-SVD can be useful in this context even when the images are quite noisy.
Moreover, as we show in Fig \ref{pm_fig_image_pair_X_wSM_FIGF}, a judiciously chosen $\Ckern$ can be used to generate principal-modes that are effective for aligning noisy images to one another.

The story for volumes is very similar; the success of the radial- and degree-SVD is linked to the decay of the spectrums of $\Ckern(\fB)$ and $\Dkern(\fB)$.
Once again, when the $\fB$ used to build these kernels is relatively smooth, then these kernels will be approximately low-rank, and the operation-count of the inner-product calculation can be reduced considerably.

\begin{figure}
\centering
\includegraphics[width=7.5in]{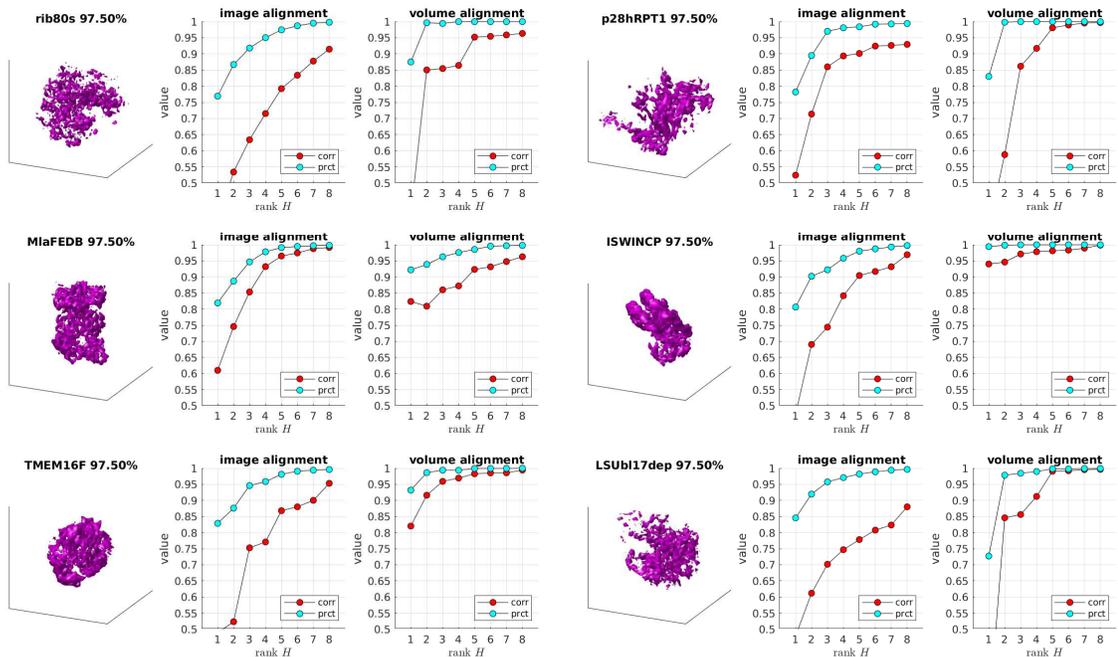}
\caption{\label{pm_vfig_collect_11}
Here we apply the same strategy to a variety of other data-sets.
In reading order we show the $97.5\%$ level-set for a published molecule (see main text for list), followed by the correlation and backwards-error for image- and volume-alignment, each as a function of $\nrank$.
}
\end{figure}

The examples in Figs \ref{pm_fig_X_2d_cluster_d0_FIGK} and \ref{pm_vfig_timing_FIGB} were constructed using the {\tt EMPIAR-10005} dataset for the TRPV1-molecule, but the same story holds for every other cryo-em data-set we have tried.
Shown in Fig \ref{pm_vfig_collect_11} is a collection of different molecules, along with the results of our image- and volume-alignment strategy.
Just as we did above, we used the first $1024$ picked particle images for image-alignment, and multiple de-novo molecular reconstructions of the molecule for volume-alignment.
The list of data-sets and molecules is:
\begin{description}
%\item[\tt EMPIAR-10005, emd-5778:] TRPV1 dataset taken on a K2 direct electron detector.
\item[\tt EMPIAR-10028, emd-2660:] Plasmodium falciparum 80S ribosome bound to the anti-protozoan drug emetine.
\item[\tt EMPIAR-10091, emd-8674:] p28-Bound Human Proteasome Regulatory Particle, State T1.
\item[\tt EMPIAR-10536, emd-22116:] MlaFEDB from E. coli in nanodisc.
\item[\tt EMPIAR-10482, emd-9718:] ISWI-NCP complex in the ADPBeF-bound state.
\item[\tt EMPIAR-10278, emd-20244:] TMEM16F in digitonin with calcium bound.
\item[\tt EMPIAR-10076, emd-8434:] L17-Depleted 50S Ribosomal Intermediates, class A.
\end{description}
Note that, in each case, a low backwards-error is obtained even when $\nrank$ is relatively small.

In this paper we have only presented the most idealized objective-functions, each of which integrates over the uniform distribution of rotations.
Our approach can immediately be generalized to other objective-functions that account for a non-uniform distribution of rotations, and/or a distribution of translations.
For example, if we consider rigid image-alignment over a distribution of rotations $\mu_{\gamma}$ on $[0,2\pi)$, and a distribution of translations $\mu_{\vd}$ on $\Real^{2}$, then the associated objective-function for a single target $\fB$ would be proportional to:
\begin{eqnarray}
%\Ccost(\vu;\fB) = \int_{0}^{2\pi} \int_{\Real^{2}} \int_{0}^{2\pi} \int_{\Real^{2}} \left\{ \sum_{q=0}^{\qmax-1} \left| \left[\transpose{\vu}\translation_{\vd}\rotation_{\gamma}\fB\right](\psi_{q}) - \left[\transpose{\vu}\translation_{\vd'}\rotation_{\gamma'}\fB\right](\psi_{q}) \right|^{2} \right\} d\mu_{\vd}(\vd) d\mu_{\gamma}(\gamma) d\mu_{\vd}(\vd') d\mu_{\gamma}(\gamma') \text{,}
\Ccost(\vu;\fB) = \iiiint \left\{ \sum_{q=0}^{\qmax-1} \left| \left[\transpose{\vu}\translation_{\vd}\rotation_{\gamma}\fB\right](\psi_{q}) - \left[\transpose{\vu}\translation_{\vd'}\rotation_{\gamma'}\fB\right](\psi_{q}) \right|^{2} \right\} d\mu_{\vd}(\vd) d\mu_{\gamma}(\gamma) d\mu_{\vd}(\vd') d\mu_{\gamma}(\gamma') \text{,}
\label{eq_Image_single_cost_with_translations}
\end{eqnarray}
where $\translation_{\vd}$ represents translation (in real-space) by vector $\vd\in\Real^{2}$.
Once again, the target-averaged objective function can be obtained simply by summing over different targets $\fB_{j}$, as in \eqref{eq_image_cost_and_kern}.
Just as before, this objective function is quadratic in $\vu$, meaning that it can be optimized by finding the principal-vectors of the associated kernel.

In the context of cryo-EM, objectives of the form \eqref{eq_Image_single_cost_with_translations} arise naturally when aligning noisy images (which are only approximately centered) to targets which are formed by projecting a reference molecule.
%Indeed, if we assume (i) that $\mu_{\gamma}$ is uniform, and (ii) that $\mu_{\vd}$ is an isotropic gaussian in $\Real^{2}$ centered at the origin with standard-deviation $\sigma_{\vd}$, and (iii) that the viewing-angles of each of the targets are uniformly distributed across $S03$, then the target-averaged kernel $\Ckern$ associated with \eqref{eq_Image_single_cost_with_translations} converges to:
Indeed, if we assume (i) that $\mu_{\vd}$ is an isotropic gaussian in $\Real^{2}$ centered at the origin with standard-deviation $\sigma_{\vd}$, and (ii) that the viewing-angles of each of the targets are uniformly distributed across $S03$, then the target-averaged kernel $\Ckern$ associated with \eqref{eq_Image_single_cost_with_translations} converges (up to a constant factor) to:
\begin{eqnarray}
\Ckern_{r,r'} & = & \left[ \sum_{l=0}^{+\infty} \sum_{m=-l}^{m=+l} \fB_{l}^{m}(k_{r})^{\dagger}\fB_{l}^{m}(k_{r'}) \right] \cdot \boldsymbol{E}^{+}(k_{r},k_{r'}) - \left[ \fB_{0}^{0}(k_{r})^{\dagger}\fB_{0}^{0}(k_{r'}) \right] \cdot \boldsymbol{E}^{-}(k_{r}) \cdot \boldsymbol{E}^{-}(k_{r'}) \text{,}
\label{eq_template_cost_with_translations}
\end{eqnarray}
as $\qmax\rightarrow\infty$, where the terms $\boldsymbol{E}^{+}$ and $\boldsymbol{E}^{-}$ denote
\begin{eqnarray}
  \boldsymbol{E}^{+}(k_{r},k_{r'}) = \exp\left(-\frac{\sigma_{\vd}^{2}}{2}\left[ k_{r}^{2} + k_{r'}^{2} \right]\right) \cdot \exp\left( k_{r}k_{r'}\sigma_{\vd}^{2} \right) \text{,}
\end{eqnarray}
and
\begin{eqnarray}
\boldsymbol{E}^{-}(k_{r}) = \tilde{k}_{r}\sqrt{\frac{\pi}{2}} \exp\left(-\tilde{k}_{r}^{2}\right) \left( {\cal I}_{-1/2}(\tilde{k}_{r}^{2}) - {\cal I}_{+1/2}(\tilde{k}_{r}^{2}) \right) \text{,}
\end{eqnarray}
where $\tilde{k}_{r} = k_{r}\sigma_{\vd}/2$, and ${\cal I}_{q}$ refers to the modified bessel function of first kind of order $q$.
We remark that the translation-factorization of \cite{RSAB20} combines naturally with the radial-SVD, further motivating an objective-function kernel of the form \eqref{eq_template_cost_with_translations}.

Finally, we comment that there are several strategies for alignment that involve coarsely sampling the landscape of inner-products, and then refining the landscape near regions that look promising \cite{Punjani2017}.
The low backwards-error of our radial-SVD indicates that this approach might be quite useful for building an approximate inner-product landscape that highlights the appropriate regions of interest.
The high correlation between our approximate inner-product landscape and the full calculation also implies that -- as an alternative to the frequency-marching proposed by \cite{Barnett2017} -- one might consider a version of `principal-mode-marching', where a low value of $\nrank$ is chosen initially to build a coarse approximation of the inner-product landscape, and then $\nrank$ is successively increased to refine the calculation.

As one can readily see, the concepts described in this paper are not limited to image- and volume-alignment, and can be applied throughout the molecular-reconstruction pipeline.
Indeed, the principal-volumes of the molecule can be reconstructed by solving a least-squares problem restricted to the same set of principal-modes that were used to align the principal-image-rings to the principal-template-rings.
We defer this discussion to future work.

\bibliographystyle{unsrt}
\bibliography{Principled_Marching_bib}

\begin{thebibliography}{10}

\bibitem{CapekPousek2003}
M.~Capek and L.~Pousek.
\newblock Biomedical volume alignment using an efficient optimization method
  and fast data resampling.
\newblock In {\em Proceedings of the 3rd IEEE International Symposium on Signal
  Processing and Information Technology (IEEE Cat. No.03EX795)}, pages
  483--486, 2003.

\bibitem{Szeliski2006}
R.~Szeliski.
\newblock Image alignment and stitching: {A} tutorial.
\newblock {\em Found. Trends Comput. Graph. Vis.}, 4(1):1--104, 2006.

\bibitem{Shatsky2009}
Maxim Shatsky, Richard~J. Hall, Steven~E. Brenner, and Robert~M. Glaeser.
\newblock {A method for the alignment of heterogeneous macromolecules from
  electron microscopy}.
\newblock {\em jsb}, 166:67--78, 2009.

\bibitem{Cheng2015}
Yifan Cheng, Nikolaus Grigorieff, Pawel~A. Penczek, and Thomas Walz.
\newblock A primer to single-particle cryo-electron microscopy.
\newblock {\em Cell}, 161:439--449, 2015.

\bibitem{Nogales2015}
E.~Nogales and S.H. Scheres.
\newblock {Cryo-EM}: a unique tool for the visualization of macromolecular
  complexity.
\newblock {\em Mol. Cell.}, 58:677--689, 2015.

\bibitem{Elmlund2015}
Dominika Elmlund and Hans Elmlund.
\newblock Cryogenic electron microscopy and single-particle analysis.
\newblock {\em Annu. Rev. Biochem.}, 84:499--517, 2015.

\bibitem{Murata2017}
K.~Murata and M.~Wolf.
\newblock Cryo-electron microscopy for structural analysis of dynamic
  biological macromolecules.
\newblock {\em Biochim. Biophys. Acta Gen. Subj.}, 1862(2):324--334, 2017.

\bibitem{Sigworth2016}
F~Sigworth.
\newblock Principles of {cryo-EM} single-partice image processing.
\newblock {\em Microscopy}, 65(1):57--67, 2016.

\bibitem{Goncharov1987}
A.~B. Goncharov.
\newblock {Methods of integral geometry and finding the relative orientation of
  identical particles arbitrarily arranged in a plane from their projections
  onto a straight line}.
\newblock {\em Dokl. Phys.}, 32:173, 1987.

\bibitem{vanHeel1987}
M.~van Heel.
\newblock {Angular reconstitution: {A} posteriori assignment of projection
  directions for 3D reconstruction}.
\newblock {\em Ultramicroscopy}, 21:111--123, 1987.

\bibitem{Goncharov1988}
A.~B. Goncharov and M.~S. Gelfand.
\newblock {Determination of mutual orientation of identical particles from
  their projections by the moments method}.
\newblock {\em Ultramicroscopy}, 25:317--328, 1988.

\bibitem{Jonic2005}
S.~Jonic, C.O. Sorzano, P.~Thevenaz, C.~El-Bez, S.~De~Carlo, and M.~Unser.
\newblock Spline-based image-to-volume registration for three-dimensional
  electron microscopy.
\newblock {\em Ultramicroscopy}, 103(4):303--317, 2005.

\bibitem{Grigorieff2007}
Nikolaus Grigorieff.
\newblock {FREALIGN: High-resolution refinement of single particle structures}.
\newblock {\em J. Struct. Biol.}, 157(1):117--125, 2007.

\bibitem{EMAN2}
G.~Tang, L.~Peng, P.R. Baldwin, D.S. Mann, W.~Jiang, I.~Rees, and S.J. Ludtke.
\newblock {EMAN2}: {A}n extensible image processing suite for electron
  microscopy.
\newblock {\em J. Struct. Biol.}, 157:38--46, 2007.

\bibitem{Yang2008}
Zh. Yang and P.A. Penczek.
\newblock {Cryo-EM} image alignment based on nonuniform fast {Fourier}
  transform.
\newblock {\em Ultramicroscopy}, 108:959--969, 2008.

\bibitem{Singer2009}
A.~Singer, R.~R. Coifman, F.~J. Sigworth, D.~W. Chester, and Y.~Shkolnisky.
\newblock Detecting consistent common lines in {cryo-EM} by voting.
\newblock {\em J. Struct. Biol.}, 169:312--322, 2009.

\bibitem{Singer2011}
A.~Singer and Y.~Shkolnisky.
\newblock Three-dimensional structure determination from common lines in
  {cryo-EM} by eigenvectors and semidefinite programming.
\newblock {\em SIAM J. Imaging Sci.}, 4(2):543--572, 2011.

\bibitem{Scheres2012}
Sjors H~W Scheres.
\newblock A {B}ayesian view on {cryo-EM} structure determination.
\newblock {\em J. Mol. Biol.}, 415:406--418, 2012.

\bibitem{Shkolnisky2012}
Y.~Shkolnisky and A.~Singer.
\newblock {Viewing direction estimation in {cryo-EM} using synchronization.}
\newblock {\em SIAM J. Imaging Sci.}, 5(3):1088--1110, 2012.

\bibitem{Lyumkis2013}
Dmitry Lyumkis, Axel~F. Brilot, Douglas~L. Theobald, and Nikolaus Grigorieff.
\newblock {Likelihood-based classification of cryo-EM images using FREALIGN}.
\newblock {\em J. Struct. Biol.}, 183(3):377--388, 2013.

\bibitem{Wang2013}
L.~Wang, A.~Singer, and Z.~Wen.
\newblock {Orientation determination from {cryo-EM} images using least
  unsquared deviations}.
\newblock {\em SIAM J. Imaging Sci.}, 6(4):2450--83, 2013.

\bibitem{Grigorieff2016}
N.~Grigorieff.
\newblock Frealign: An exploratory tool for single-particle {cryo-EM}.
\newblock In R~A Crowther, editor, {\em The Resolution Revolution: Recent
  Advances In {cryoEM}}, volume 579 of {\em Methods Enzymol.}, pages 191--226.
  Academic Press, 2016.

\bibitem{Punjani2017}
A.~Punjani, J.L. Rubinstein, D.J. Fleet, and M.A. Brubaker.
\newblock {cryoSPARC}: algorithms for rapid unsupervised {cryo-EM} structure
  determination.
\newblock {\em Nat. Methods}, 14:290--296, 2017.

\bibitem{PunjaniBrubaker2017}
A.~Punjani, M.A. Brubaker, and D.J. Fleet.
\newblock Building proteins in a day: efficient {3D} molecular structure
  estimation with electron cryomicroscopy.
\newblock {\em IEEE Trans. Pattern Anal. Mach. Intell.}, 39(4):706--718, 2017.

\bibitem{RSAB20}
Aaditya Rangan, Marina Spivak, Joakim And{\'e}n, and Alex Barnett.
\newblock Factorization of the translation kernel for fast rigid image
  alignment.
\newblock {\em Inverse Problems}, 36(2), 2020.

\bibitem{Joyeux2002}
Laurent Joyeux and Pawel~A. Penczek.
\newblock {Efficiency of 2D alignment methods}.
\newblock {\em Ultramicroscopy}, 92(2):33--46, 2002.

\bibitem{Barnett2017}
A.~Barnett, L.~Greengard, A.~Pataki, and M.~Spivak.
\newblock Rapid solution of the {cryo-EM} reconstruction problem by frequency
  marching.
\newblock {\em SIAM J. Imaging Sci.}, 10(3):1170--1195, 2017.

\bibitem{Kostelec03fftson}
Peter~J. Kostelec and Daniel~N. Rockmore.
\newblock Ffts on the rotation group.
\newblock In {\em Santa Fe Institute Working Papers Series Paper}, pages
  03--11, 2003.

\bibitem{bracewell}
R~Bracewell.
\newblock {\em The {F}ourier Transform and Its Applications}.
\newblock McGraw-Hill, 3rd edition, 1999.

\bibitem{Zhao2014}
Zh. Zhao and A.~Singer.
\newblock Rotationally invariant image representation for viewing direction
  classification in {cryo-EM}.
\newblock {\em J. Struct. Biol.}, 186:153--166, 2014.

\bibitem{Zhao2016}
Zh. Zhao, Y.~Shkolnisky, and A.~Singer.
\newblock Fast steerable principal component analysis.
\newblock {\em IEEE Trans. Comput. Imaging}, 2(1):1--12, 2016.

\bibitem{finufft}
A~H Barnett, J~F Magland, and L~af~Klinteberg.
\newblock A parallel non-uniform fast {F}ourier transform library based on an
  ``exponential of semicircle'' kernel, 2019.

\bibitem{Liao2013}
M.~Liao, E.~Cao, D.~Julius, and Y.~Cheng.
\newblock {Structure of the TRPV1 ion channel determined by electron
  cryo-microscopy}.
\newblock {\em Nature}, 504:107--12, 2013.

\bibitem{Sigworth1998}
F~J Sigworth.
\newblock {A maximum-likelihood approach to single-particle image refinement}.
\newblock {\em J. Struct. Biol.}, 122(3):328--39, 1998.

\bibitem{Scheres2009}
S~H~W Scheres, M~Valle, P~Grob, E~Nogales, and J.-M. Carazo.
\newblock Maximum likelihood refinement of electron microscopy data with
  normalization errors.
\newblock {\em J. Struct. Biol.}, 166(2):234--240, 2009.

\bibitem{Sigworth2010}
F.~J. Sigworth, Doerschuk P.C., J.-M. Carazo, and S.H.W. Scheres.
\newblock {An introduction to maximum-likelihood methods in Cryo-EM}.
\newblock In {\em Methods in Enzymology. Cryo-EM, Part B: 3D reconstruction},
  pages 263--294. Academic Press., 2010.

\bibitem{Baxter2009}
W.~Baxter, R.~Grassucci, Haixiao Gao, and J.~Frank.
\newblock Determination of signal-to-noise ratios and spectral snrs in cryo-em
  low-dose imaging of molecules.
\newblock {\em Journal of structural biology}, 166 2:126--32, 2009.

\bibitem{Scheres2012b}
Sjors H~W Scheres.
\newblock {RELION: Implementation of a Bayesian approach to cryo-EM structure
  determination}.
\newblock {\em J. Struct. Biol.}, 180(3):519--530, 2012.

\bibitem{Kimanius2016}
Dari Kimanius, Bjorn~O Forsberg, Sjors~HW Scheres, and Erik Lindehl.
\newblock {Accelerated cryo-EM structure determination with parallelisation
  using GPUs in RELION-2}.
\newblock {\em eLife}, 5:e18722, 2016.

\bibitem{Bell2016}
James~M. Bell, Muyuan Chen, Philip~R. Baldwin, and Steven~J. Ludtke.
\newblock {High resolution single particle refinement in EMAN2.1}.
\newblock {\em Methods}, 100:25--34, 2016.

\bibitem{BL81}
L.~C. Biedenharn and J.~D. Louck.
\newblock {\em Angular Momentum in Quantum Physics}.
\newblock Addison-Wesley, 1981.

\bibitem{Henderson2012}
Richard Henderson, Andrej Sali, Matthew~L. Baker, Bridget Carragher, Batsal
  Devkota, Kenneth~H. Downing, Edward~H. Egelman, Zukang Feng, Joachim Frank,
  Nikolaus Grigorieff, Wen Jiang, Steven~J. Ludtke, Ohad Medalia, Pawel~A.
  Penczek, Peter~B. Rosenthal, Michael~G. Rossmann, Michael~F. Schmid,
  Gunnar~F. Schr{\"{o}}der, Alasdair~C. Steven, David~L. Stokes, John~D.
  Westbrook, Willy Wriggers, Huanwang Yang, Jasmine Young, Helen~M. Berman, Wah
  Chiu, Gerard~J. Kleywegt, and Catherine~L. Lawson.
\newblock {Outcome of the first electron microscopy validation task force
  meeting}.
\newblock {\em Structure}, 20(2):205--214, 2012.

\bibitem{Penczek2014}
Pawel~A. Penczek.
\newblock {Ab initio cryo-EM structure determination as a validation problem}.
\newblock 2014.

\bibitem{Heymann2015}
J~Bernard Heymann.
\newblock {Validation of 3D EM Reconstructions: The Phantom in the Noise.}
\newblock {\em AIMS Biophysics}, 2:21--35, 2015.

\bibitem{Rosenthal2015}
Peter~B. Rosenthal and John~L. Rubinstein.
\newblock {Validating maps from single particle electron cryomicroscopy}.
\newblock {\em Current Opinion in Structural Biology}, 34:135--144, 2015.

\end{thebibliography}

\end{document}